\documentclass[10pt,a4paper]{article}

\usepackage{tikz}
\usetikzlibrary{patterns,arrows}
\usepackage{amsmath,amssymb,bbm}
\usepackage{theorem}

\newtheorem{prop}{Proposition}[section]
\newtheorem{lemma}[prop]{Lemma}
\newtheorem{theo}[prop]{Theorem}
\newtheorem{coroll}[prop]{Corollary}
{\theorembodyfont{\rmfamily}
 \newtheorem{remark}[prop]{Remark}
 \newtheorem{definition}[prop]{Definition}
 
}
\newcommand{\theoremend}{\mbox{}\hfill$\Box$}
\newcommand{\theoremendup}{\vspace{-1.2em}\theoremend}
\newcommand{\theoremendpar}{\theoremend\\[1.5ex]}
\newcommand{\theoremendskip}{\theoremend\vspace{1.5ex}}
\newcommand{\theoremendupskip}{\theoremendup\vspace{1.5ex}}
\newcommand{\definend}{\mbox{}\hfill$\lrcorner$}
\newcommand{\definendup}{\vspace{-1.2em}\definend}

\newcommand{\pmat}[1]{\begin{pmatrix}#1\end{pmatrix}}

\DeclareMathOperator{\Real}{Re}

\DeclareMathOperator{\dist}{dist}

\newcommand{\bbN}{\mathbbm{N}}
\newcommand{\bbZ}{\mathbbm{Z}}
\newcommand{\bbR}{\mathbbm{R}}
\newcommand{\bbC}{\mathbbm{C}}
\newcommand{\mdef}{\mathcal{D}}
\newcommand{\range}{\mathcal{R}}
\newcommand{\eps}{\varepsilon}
\newcommand{\rsub}{\mathcal{L}}

\begin{document}

\title{Riesz Bases for $p$-Subordinate Perturbations of Normal 
  Operators
}
\author{Christian Wyss\footnote{University of Bern, Mathematical Institute,
    Sidlerstrasse 5, 3012 Bern, Switzerland, \texttt{cwyss@math.unibe.ch}}}
\date{\today}
\maketitle

\begin{center}
  \parbox{11cm}{\small \textbf{Abstract.}
    For $p$-subordinate
    perturbations of unbounded normal operators,
    the change of the spectrum is studied
    and spectral criteria for the existence of a
    Riesz basis with parentheses of root vectors are established.
    A Riesz basis without parentheses is obtained
    under an additional a priori assumption on the spectrum of 
    the perturbed operator.
    The results are applied to two classes of
    block operator matrices.} \\[1.5ex]
  \parbox{11cm}{\small\textbf{Keywords.} Perturbation theory,
    subordinate perturbation, Riesz basis, eigenvector expansion, 
    spectrum.} \\[1.5ex]
  \parbox{11cm}{\small\textbf{Mathematics Subject Classification.}  
    Primary 47A55; Secondary 47A10, 47A70.}
\end{center}
\vspace{0em}

\section{Introduction}

Since for non-normal operators there is no analogue of the spectral theorem,
the existence of a Riesz basis (possibly with parentheses)
of root vectors is an important property:
it allows e.g.\ the construction of non-trivial invariant subspaces
and yields spectral criteria related to semigroup generation.
For a class of non-normal perturbations of normal operators
we establish different conditions in terms of the spectrum which imply
the existence of such Riesz bases.
Our assumptions on the multiplicities of the eigenvalues
are weaker than in classical perturbation theorems.

We consider an unbounded operator $T=G+S$ on a Hilbert space where 
$G$ is normal with compact resolvent and $S$ is \emph{$p$-subordinate} to $G$,
i.e.
\[\|Su\|\leq b\|u\|^{1-p}\|Gu\|^p \quad\text{for all}\quad u\in\mdef(G)\]
where $p\in[0,1[$ and $b\geq0$.
In Theorem~\ref{theo_psubpert_invl2decomp} we prove that $T$ admits
a Riesz basis with parentheses of root vectors if the eigenvalue
multiplicities of $G$ satisfy a certain asymptotic growth condition.
This growth condition is weaker than the one in a similar
result by Markus and Matsaev 
\cite{markus-matsaev}, \cite[Theorem~6.12]{markus}.

In Theorem~\ref{theo_psubpert_rootl2decomp} we obtain a Riesz basis with 
parentheses under a spectral condition of different type:
we impose no restriction on the multiplicities and instead assume that
the eigenvalues of $G$ lie on sufficiently separated line segments,
see Figure~\ref{figure_theopsub}.
If we know a priori that the eigenvalues of the perturbed operator $T$
are uniformly separated, then Theorem~\ref{theo_psubpert_rootl2decomp}
even yields a Riesz basis without parentheses. An example for such a
situation may be found in Theorem~\ref{theo_uposblkop}.
In contrast to our result, classical perturbation theorems for Riesz
bases without parentheses such as
Kato \cite[Theorem~V.4.15a]{kato},
Dunford and Schwartz \cite[Theorem~XIX.2.7]{dunford-schwartz3} and
Clark~\cite{clark}
require that almost all eigenvalues of $G$ are simple.

Apart from the above mentioned theorems,
a wide range of existence results for Riesz bases of root vectors
may be found in the literature,
both for abstract operator settings and for concrete applications.
Dissipative operators, for example, were considered by several authors;
references and some results may be found in \cite{gohberg-krein}.
For generators of $C_0$-semigroups, a Riesz basis of
eigenvectors implies the so-called spectrum determined growth assumption,
see \cite[Theorem~2.3.5]{curtain-zwart}.
Zwart~\cite{zwart} obtained Riesz bases for generators of $C_0$-groups,
while Xu and Yung~\cite{xu-yung} constructed
Riesz bases with parentheses for semigroup generators.
Riesz basis properties of root vectors are also investigated for
operator pencils, see e.g.\ \cite{azizov-dijksma-sukhocheva,vizitei-markus}.
Pencils coming from concrete physical problems were studied in
\cite{adamjan-pivovarchik-tretter,marletta-shkalikov-tretter}.

Finally there are simple examples of non-normal operators
whose eigenvectors are complete but do not form a Riesz basis,
see e.g.\ \cite{davies99}.

In this paper, we follow ideas due to 
Markus and Matsaev~\cite[Chapter~1]{markus} to prove the existence of 
Riesz bases of root vectors.
In Section~\ref{sec_complete} we start by deriving a completeness theorem 
for the system of root vectors of an operator with compact resolvent.
Unlike the classical Keldysh theorem on
completeness \cite{keldysh}, \cite[\S4]{markus}, where the resolvent
belongs to a von~Neumann-Schatten class, we assume here that
it is uniformly bounded on an appropriate sequence of curves.

In Section~\ref{sec_rbasis} we then recall the notion of a Riesz basis
consisting of subspaces and provide a sufficient condition for 
its existence in terms of projections.
Although a Riesz basis consisting of (finite-dimensional) subspaces
is equivalent to a Riesz basis with parentheses,
we use the basis of subspaces notion in the formulation of our theorems,
since it is more convenient.

In Section~\ref{sec_psubspec} we study in detail the change of the spectrum 
of a normal operator under a $p$-subordinate perturbation. 
The basic observation here is that if the spectrum of $G$ lies on rays
from the origin, then the spectrum of $T$ lies inside parabolas around
these rays.
Based on the localisation of the spectrum, several estimates for 
Riesz projections of $T$ are obtained in Section~\ref{sec_rprojest}.

In Section~\ref{sec_psubpert} we derive our main existence results
for Riesz bases of root vectors.
In fact, these results also hold in the more general setting where $G$ 
is a possibly non-normal operator with compact resolvent, a Riesz basis of 
root vectors and an appropriate spectrum,
see Proposition~\ref{prop_psubpert_revert} and
Remark~\ref{rem_psubpert_revert}.

In Section~\ref{sec_psubpert_appli} we finally apply our theory to
diagonally dominant block operator matrices.
Riesz bases of root vectors were obtained by
Jacob, Trunk and Winklmeier~\cite{jacob-trunk-winklmeier} for
operator matrices associated with damped vibrations of thin beams
and by Kuiper and Zwart~\cite{kuiper-zwart} for a class of Hamiltonian
operator matrices from control theory.
In Theorem~\ref{theo_psubblkop} we consider operator matrices whose entries
may all be unbounded, whereas in
\cite{jacob-trunk-winklmeier,kuiper-zwart} some of the entries
were always bounded.
Theorem~\ref{theo_uposblkop} applies to a class of Hamiltonians which 
is different from the one in \cite{kuiper-zwart}.
While the eigenvalues of the diagonal part of the Hamiltonian in
\cite{kuiper-zwart} are simple, we consider the case of double eigenvalues.

\section{Completeness of the system of root subspaces}
\label{sec_complete}

We derive a completeness theorem for the system of root subspaces
of an operator with compact resolvent, which applies to a
different situation than the classical theorem of Keldysh
\cite{keldysh}, \cite[\S4]{markus}.

Let $T$ be an operator on a Banach space with a compact isolated part
$\sigma\subset\sigma(T)$ of the spectrum.
Let $\Gamma$ be a simply closed, positively oriented integration contour
with $\sigma$ in its interior and $\sigma(T)\setminus\sigma$ in its exterior.
Then
\begin{equation}\label{rieszproj_def}
  P=\frac{i}{2\pi}\int_\Gamma(T-z)^{-1}\,dz
\end{equation}
defines a projection such that
$\range(P)$ and $\ker P$ are
$T$-invariant, $\range(P)\subset\mdef(T)$, and
\[\sigma(T|_{\range(P)})=\sigma,\qquad
  \sigma(T|_{\ker P})=\sigma(T)\setminus\sigma.\]
$P$ does not depend on the particular choice of $\Gamma$
and is called the \emph{Riesz projection} associated with 
$\sigma$ (or $\Gamma$); $\range(P)$ is the corresponding 
\emph{spectral subspace}.
For a proof see
\cite[Theorem~XV.2.1]{gohberg-goldberg-kaashoek} or
\cite[Theorem~III.6.17]{kato}.

For an eigenvalue $\lambda$ of
$T$ we call
\[\rsub(\lambda)=\bigcup_{k\in\bbN}\ker(T-\lambda)^k\]
the \emph{root subspace} of $T$ corresponding to $\lambda$;
the non-zero elements of $\rsub(\lambda)$ are the \emph{root vectors}. 
A sequence of root vectors
$x_1,\dots,x_n\in\rsub(\lambda)$ is called a \emph{Jordan chain}
if $(T-\lambda)x_k=x_{k-1}$ for $k\geq 2$ and $(T-\lambda)x_1=0$.
In the case that $T$ has a compact resolvent, its spectrum consists
of isolated eigenvalues only; so for every eigenvalue $\lambda$ there
is the associated Riesz projection $P_\lambda$, which satisfies
$\range(P_\lambda)=\rsub(\lambda)$.

Recall that for an operator $T$ with compact resolvent on a
Hilbert space its adjoint $T^*$ also has a compact resolvent.
\begin{lemma}
  Let\, $T$ be an operator with compact resolvent
  on a Hilbert space and
  $M$ the subspace generated by all root subspaces of\, $T$, i.e.,
  the set of all finite linear combinations of root vectors of\, $T$.
  If\, $P$ is the Riesz projection of\, $T^*$
  corresponding to an eigenvalue $\lambda\in\sigma(T^*)$, then
  $M^\perp\subset\ker P$.
  Moreover, $M^\perp$ is $T^*$-invariant and
  $(T^*-z)^{-1}$-invariant for every $z\in\varrho(T^*)$; in particular
  $\varrho(T^*)\subset\varrho(T^*|_{M^\perp})$.
\end{lemma}
\textit{Proof.} 
We have $\lambda\in\sigma(T^*)$ if and only if 
$\overline{\lambda}\in\sigma(T)$.
Observe that if $P$ is the Riesz projection of $T^*$ corresponding
to $\lambda$, then $P^*$ is the Riesz projection of $T$ corresponding
to $\overline{\lambda}$.
Since $\range(P^*)\subset M$ we find
\(M^\perp\subset\range(P^*)^\perp=\ker P\).
Now let $v\in M$ and $z\in\varrho(T^*)$. Then 
$Tv,(T-\bar{z})^{-1}v\in M$ and we have
\begin{align*}
  u\in M^\perp\cap\mdef(T^*)\quad&\Rightarrow\quad(T^*u|v)=(u|Tv)=0,\\
  u\in M^\perp\quad&\Rightarrow\quad
  \bigl((T^*-z)^{-1}u\big|v\bigr)=\bigl(u\big|(T-\bar{z})^{-1}v\bigr)=0.
\end{align*}
Therefore $M^\perp$ is $T^*$- and $(T^*-z)^{-1}$-invariant, and
this in turn implies the inclusion
$\varrho(T^*)\subset\varrho(T^*|_{M^\perp})$.
\theoremend

\begin{coroll}
  Let\, $T$ and $M$ be as above. Then $\varrho(T^*|_{M^\perp})=\bbC$.
\end{coroll}
\textit{Proof.} Since $T$ has a compact resolvent, the same holds for
$T^*$ and $T^*|_{M^\perp}$. Consequently if 
$\lambda\in\sigma(T^*|_{M^\perp})$, then $\lambda$ is an eigenvalue
of $T^*|_{M^\perp}$, i.e., $T^*u=\lambda u$ for some 
$u\in M^\perp\setminus\{0\}$. In particular $\lambda$ is an
eigenvalue of $T^*$ and we have $u\in\range(P)$ where $P$ is the
Riesz projection of $T^*$ corresponding to $\lambda$.
Now the previous lemma implies $u\in M^\perp\subset\ker P$
and hence $u=0$, which is a contradiction. Therefore
$\sigma(T^*|_{M^\perp})=\varnothing$.
\theoremend

\begin{theo}\label{theo_complrsub}
  Let\, $T$ be an operator with compact resolvent on a Hilbert space $H$
  with scalar product $(\cdot|\cdot)$.
  If there exists a sequence of bounded regions $(U_k)_{k\in\bbN}$ such
  that $\bbC=\bigcup_{k\in\bbN}U_k$, $\partial U_k\subset\varrho(T)$
  for all $k$, and there is a constant $C\geq 0$ with
  \[\|(T-z)^{-1}\|\leq C \quad\text{for}\quad z\in\partial U_k,\,k\in\bbN,\]
  then the system of root subspaces of\, $T$ is complete.\footnote{%
      A system of subspaces in $H$ is called \emph{complete} if the subspace
      generated by the system is dense in $H$.}
\end{theo}
\textit{Proof.} Let $M$ be as before.
For $u,v\in M^\perp$ we consider the holomorphic function defined by 
\[f(z)=\bigl((T^*|_{M^\perp}-z)^{-1}u|v\bigr).\]
From the previous corollary
we know that its domain of definition is $\bbC$. Since
\[\|(T^*|_{M^\perp}-z)^{-1}\|\leq\|(T^*-z)^{-1}\|=\|(T-\bar{z})^{-1}\|
  \quad\text{for}\quad \bar{z}\in\varrho(T),\]
we see that $|f(z)|\leq C\|u\|\|v\|$ holds for
$\bar{z}\in\partial U_k$. Using the maximum principle,
we find that $|f(z)|\leq C\|u\|\|v\|$
for every $z\in\bbC$; by Liouville's theorem
$f$ is constant. Since $u$ and $v$ have been arbitrary, the mapping
$z\mapsto(T^*|_{M^\perp}-z)^{-1}$ is also constant. For $u\in M^\perp$
we obtain
\begin{align*}
  (T^*|_{M^\perp})^{-1}u=(T^*|_{M^\perp}-I)^{-1}u \quad&\Rightarrow\quad
  (T^*|_{M^\perp}-I)(T^*|_{M^\perp})^{-1}u=u\\
  \Rightarrow\quad (T^*|_{M^\perp})^{-1}u=0 \quad&\Rightarrow\quad u=0\,.
\end{align*}
Hence $M^\perp=\{0\}$, i.e., $M\subset H$ is dense.
\theoremend

\begin{coroll}\label{coroll_complete_rootspace}
  Let\, $T$ be an operator with compact resolvent on a Hilbert space.
  Suppose that almost all eigenvalues of\, $T$ 
  lie in a finite number of pairwise disjoint sectors
  \[\Omega_j=\bigl\{z\in\bbC\,\big|\,|\arg z-\theta_j|<\psi_j\bigr\}
    \quad\text{with}\quad
    0<\psi_j\leq\frac{\pi}{4},\quad j=1,\ldots,n.\]
  If there are constants $C,r_0\geq0$ such that
  \[\|(T-z)^{-1}\|\leq C\quad\text{for}\quad z\not\in\Omega_1\cup\ldots\cup
    \Omega_n,\,|z|\geq r_0\]
  and for each sector $\Omega_j$ there is a sequence $(x_k)_{k\in\bbN}$
  with $x_k\to\infty$ and
  \[\|(T-z)^{-1}\|\leq C\quad\text{for}\quad z\in\Omega_j,\,
    \Real(e^{-i\theta_j}z)=x_k,\,k\in\bbN,\]
  then the system of root subspaces of\, $T$ is complete.
  \theoremend
\end{coroll}

\section{Riesz bases of subspaces}
\label{sec_rbasis}

We recall the closely related concepts of Riesz bases,
Riesz bases with parentheses, and Riesz bases of subspaces, see
\cite[\S 1]{vizitei-markus}, \cite[Chapter VI]{gohberg-krein},
\cite[\S15]{singer2} and \cite[\S 2]{wyss} for more details.

\begin{definition}
  Let $H$ be a separable Hilbert space.
  \begin{itemize}
  \item[(i)] A sequence $(v_k)_{k\in\bbN}$ in $H$ is called a 
    \emph{Riesz basis}
    of $H$ if there is an isomorphism $J:H\to H$ such that
    $(Jv_k)_{k\in\bbN}$ is an orthonormal basis of $H$.
  \item[(ii)]
    A sequence of closed subspaces $(V_k)_{k\in\bbN}$ of $H$ is called a
    \emph{Riesz basis of subspaces} of $H$ if there is an isomorphism
    $J:H\to H$ such that $(J(V_k))_{k\in\bbN}$ is a
    complete system
    of pairwise orthogonal subspaces.
  \end{itemize}
  
  \definendup
\end{definition}
Other notions for Riesz bases of subspaces are ``unconditional basis
of subspaces'' \cite{vizitei-markus} or ``$l^2$-decomposition''
\cite{singer2}.

The sequence $(v_k)_{k\in\bbN}$ is a Riesz basis if and only if
$\inf\|v_k\|>0$, $\sup\|v_k\|<\infty$, and
every $x\in H$ has a unique representation
\[x=\sum_{k=0}^\infty\alpha_k v_k,\qquad \alpha_k\in\bbC,\]
where the convergence of the series is unconditional.
There is a similar characterisation for Riesz bases of subspaces,
see \cite[\S 2.2]{wyss} and \cite[\S VI.5]{gohberg-krein} for a proof:

\begin{prop}\label{prop_rbasissub-equiv}
  For a sequence $(V_k)_{k\in \bbN}$ of closed 
  subspaces of $H$ the following assertions are equivalent:
  \begin{itemize}
  \item[(i)] $(V_k)_{k\in\bbN}$ is a Riesz basis of subspaces for $H$.
  \item[(ii)] 
    The sequence $(V_k)_{k\in\bbN}$ is complete and
    there exists $c\geq 1$ such that
    \begin{equation}\label{l2decomp_def_ineq}
      c^{-1}\sum_{k\in F}\|x_k\|^2 \leq \Big\|\sum_{k\in F}x_k\Big\|^2
      \leq c\sum_{k\in F}\|x_k\|^2
    \end{equation}
    for all finite subsets $F\subset\bbN$ and $x_k\in V_k$.
  \item[(iii)] Every $x\in H$ has a unique representation
    $x=\sum_{k=0}^\infty x_k$ with $x_k\in V_k$,
    where the convergence of the series is unconditional.
  \end{itemize}

  \theoremendup
\end{prop}
To refer to the constant in (\ref{l2decomp_def_ineq}),
we shall also speak of a \emph{Riesz basis of subspaces with constant $c$}.

A sequence $(v_k)_{k\in\bbN}$ in a Hilbert space $H$ is called a
\emph{Riesz basis with parentheses} if there exists
a Riesz basis of subspaces $(V_k)_{k\in\bbN}$ of $H$ and
a subsequence $(n_k)_k$ of $\bbN$ with $n_0=0$ such that
$(v_{n_k},\dots,v_{n_{k+1}-1})$ is a basis of $V_k$.
In this case every $x\in H$ has a unique representation
\[x=\sum_{k=0}^\infty \Biggl(\sum_{j=n_k}^{n_{k+1}-1}\alpha_j v_j\Biggr),
  \qquad \alpha_j\in\bbC,\]
where the series over $k$ converges unconditionally.

The definition of a Riesz basis of subspaces generalises naturally
to a family of closed subspaces $(V_k)_{k\in\Lambda}$, where the index set
$\Lambda$ is either finite or countably infinite;
Proposition~\ref{prop_rbasissub-equiv} continuous to hold in this context.
In particular, a finite family $(V_1,\dots,V_n)$ of closed subspaces is a
Riesz basis of $H$ if and only if the subspaces form a direct sum
$H=V_1\oplus\dots\oplus V_n$.
Despite this equivalence, the Riesz basis notion is convenient even
for finite families
to specify the constant $c$ in \eqref{l2decomp_def_ineq}.
An example is the next lemma, which is used in the proof of 
Theorem~\ref{theo_psubpert_rootl2decomp}
to show that the root subspaces of an operator form a Riesz basis.

\begin{lemma}\label{lem_join_l2decomp}
  Let\, $(W_k)_{k\in \Lambda}$ be a Riesz basis of subspaces of\, $H$
  with constant\, $c_0$.
  Let\, $(V_{kj})_{j\in J_k}$ be 
  Riesz bases of subspaces of $W_k$ for all\, $k\in\Lambda$ with
  common constant\, $c_1$.
  Then the family
  $(V_{kj})_{k\in \Lambda,\,j\in J_k}$ is a
  Riesz basis of subspaces of\, $H$ with constant\, $c_0c_1$.
\end{lemma}
\textit{Proof.} Since $(W_k)_{k\in\Lambda}$ is complete in $H$ 
and $(V_{kj})_{j\in J_k}$ is complete 
in $W_k$ for every $k\in \Lambda$, the family
$(V_{kj})_{k\in \Lambda,j\in J_k}$ is complete in $H$. 
Consider $F\subset \Lambda$ finite, 
$F_k\subset J_k$ finite for each $k\in F$,
and $x_{kj}\in V_{kj}$. Then, using \eqref{l2decomp_def_ineq}, we obtain
\[\bigg\|\sum_{\begin{subarray}{l}k\in F\\j\in F_k
  \end{subarray}}
  x_{kj}\bigg\|^2 \leq c_0\sum_{k\in F}\bigg\|
  \sum_{j\in F_k}x_{kj}\bigg\|^2
  \leq c_0\sum_{k\in F}c_1\sum_{j\in F_k}
  \|x_{kj}\|^2
  = c_0c_1\sum_{\begin{subarray}{l}k\in F\\j\in F_k
  \end{subarray}}
  \|x_{kj}\|^2\]
and similarly
\(\|\sum_{k\in F,j\in F_k}x_{kj}\|^2\geq
  c_0^{-1}c_1^{-1}\sum_{k\in F,j\in F_k}
  \|x_{kj}\|^2\).
\theoremendpar
Note that the existence of the common constant $c_1$ is guaranteed
if only finitely many $J_k$ consist of more than one element.

Our next aim is to derive a sufficient condition
for a sequence of projections to  generate a Riesz basis of subspaces.
\begin{lemma}\label{lem_rearrange_sign}
  Let $(x_k)_{k\in\bbN}$ be a sequence in a Banach space. If
  there exists $C\geq 0$ such that for every reordering 
  \(\phi:\bbN\xrightarrow{\mathrm{bij}}\bbN\) and every $n\in\bbN$
  we have \(\|\sum_{k=0}^nx_{\phi(k)}\|\leq C\), then
  \[\sup_{n\in\bbN,\eps_k=\pm 1}\bigg\|\sum_{k=0}^n
    \eps_k x_k\bigg\|\leq 2C.\]
\end{lemma}
\textit{Proof.} Let $\eps_0,\ldots,\eps_n\in\{-1,1\}$
and consider reorderings $\phi_1$ and $\phi_2$ that move all $+1$ and
all $-1$ in the sequence $(\eps_0,\ldots,\eps_n)$, respectively, 
to its beginning.
Then, with $n_1$, $n_2$ appropriate, we obtain
\[\bigg\|\sum_{k=0}^n\eps_k x_k\bigg\|
  \leq\bigg\|\sum_{\substack{k=0\\\eps_k=+1}}^nx_k\bigg\|
  +\bigg\|\sum_{\substack{k=0\\\eps_k=-1}}^nx_k\bigg\|
  =\bigg\|\sum_{k=0}^{n_1}x_{\phi_1(k)}\bigg\|
  +\bigg\|\sum_{k=0}^{n_2}x_{\phi_2(k)}\bigg\| \leq 2C.\]

\theoremendup

\begin{lemma}
  Let\, $H$ be a Hilbert space, $x_0,\ldots,x_n\in H$, and 
  \[E=\bigl\{(\eps_0,\ldots,\eps_n)\,\big|\,\eps_k=\pm 1\bigr\}.\]
  Then
  \[2^{n+1}\sum_{k=0}^n\|x_k\|^2=\sum_{\eps\in E}
    \|\eps_0x_0+\cdots+\eps_nx_n\|^2.\]
\end{lemma}
\textit{Proof.} We use induction on $n$. The statement is true for 
the case $n=0$ since
$2\|x_0\|^2=\|x_0\|^2+\|-x_0\|^2$.
Now suppose the statement holds for some $n\geq 0$; let 
\[\widetilde{E}
  =\bigl\{(\eps_0,\ldots,\eps_{n+1})\,\big|\,\eps_k=\pm 1\bigr\}\]
and write 
$x_\eps=\eps_0x_0+\cdots+\eps_nx_n$. Then
\begin{align*}
  &\sum_{\eps\in\widetilde{E}}\|\eps_0x_0+\cdots+\eps_{n+1}x_{n+1}\|^2
  =\sum_{\eps\in E}\left(\|x_\eps+x_{n+1}\|^2+\|x_\eps-x_{n+1}\|^2\right)\\
  =\,\,&\sum_{\eps\in E}\left(2\|x_\eps\|^2+2\|x_{n+1}\|^2\right)
  =2\sum_{\eps\in E}\|x_\eps\|^2+2\cdot 2^{n+1}\|x_{n+1}\|^2\\
  =\,\,&2^{n+2}\left(\sum_{k=0}^n\|x_k\|^2+\|x_{n+1}\|^2\right).
\end{align*}

\theoremendup

\begin{lemma}\label{lem_projection_estim}
  Let\, $P_0,\ldots,P_n$ be projections in a Hilbert space $H$ with
  $P_jP_k=0$ for $j\neq k$. Then
  \[C^{-2}\sum_{k=0}^n\|P_kx\|^2\leq
    \bigg\|\sum_{k=0}^nP_kx\bigg\|^2\leq C^2\sum_{k=0}^n\|P_kx\|^2
    \quad\text{for all}\quad x\in H\]
  where $C=\max\bigl\{\|\sum_{k=0}^n\eps_kP_k\|\:\big|\:\eps_k=\pm 1\bigr\}$.
\end{lemma}
\textit{Proof.} We write $x_k=P_kx$ and use the last lemma considering 
that $\eps\in E$ for which $\|\eps_0x_0+\cdots+\eps_nx_n\|$ becomes maximal.
Then we obtain
\[\sum_{k=0}^n\|P_kx\|^2 \leq\|\eps_0x_0+\cdots+\eps_nx_n\|^2
  =\bigg\|\biggl(\sum_{k=0}^n\eps_kP_k\biggr)\biggl(\sum_{k=0}^nx_k\biggr)
  \bigg\|^2 \leq C^2\bigg\|\sum_{k=0}^nP_kx\bigg\|^2.\]
On the other hand, if we choose $\eps\in E$ such that 
$\|\eps_0x_0+\cdots+\eps_nx_n\|$ is minimal, we find
\begin{align*}
  \bigg\|\sum_{k=0}^nP_kx\bigg\|^2 &=\bigg\|\biggl(\sum_{k=0}^n\eps_kP_k\biggr)
  \biggl(\sum_{k=0}^n\eps_kx_k\biggr)\bigg\|^2\\
  &\leq C^2\,\|\eps_0x_0+\cdots+\eps_nx_n\|^2\leq C^2\sum_{k=0}^n\|P_kx\|^2.
\end{align*}

\theoremendupskip

The following statement is a slight modification of a result\footnote{
  Under the weaker assumption
  $\sum_{k=0}^\infty|(P_kx|y)|<\infty$ for all $x,y\in H$,
  the existence of the Riesz basis of subspaces is proved, but
  without obtaining an estimate for the constant $c$.}
in the book of Markus \cite[Lemma 6.2]{markus}.
\begin{prop}\label{prop_proj_l2decomp}
  Let\, $H$ be a Hilbert space and
  $(P_k)_{k\in\bbN}$ a sequence of projections in
  $H$ satisfying $P_jP_k=0$ for $j\neq k$. Suppose that 
  the family \((\range(P_k))_{k\in\bbN}\) is complete in $H$
  and that
  \begin{equation}\label{prop_proj_l2decomp-ie}
    \sum_{k=0}^\infty|(P_kx|y)|\leq C\|x\|\|y\|
    \quad\text{for all}\quad x,y\in H
  \end{equation}
  with some constant\, $C\geq 0$. Then $(\range(P_k))_{k\in\bbN}$
  is a Riesz basis of subspaces of $H$ with constant\, $c=4C^2$.
\end{prop}
\textit{Proof.} From
\[\Big|\Big(\sum_{k=0}^nP_kx\Big|y\Big)\Big| \leq
  \sum_{k=0}^n|(P_kx|y)|\leq C\|x\|\|y\|\]
we conclude that $\|\sum_{k=0}^nP_k\|\leq C$ for all $n\in\bbN$. This 
assertion remains valid after an arbitrary rearrangement of the 
sequence $(P_k)_{k\in\bbN}$ since \eqref{prop_proj_l2decomp-ie} still
holds for the rearranged sequence. An application of 
Lemmas~\ref{lem_rearrange_sign}, \ref{lem_projection_estim} and
Proposition~\ref{prop_rbasissub-equiv}
now completes the proof.
\theoremendskip

We end this section with a remark on the connection between 
Riesz bases of finite-dimensional invariant subspaces of an operator and
Riesz bases with parentheses of root vectors, see also 
\cite[\S2.3]{wyss}.
\begin{remark}\label{rem_rbasisinv-vs-parenth}
  Let $T$ be an operator on a Hilbert space. 
  Since every finite-dimen\-sional $T$-invariant subspace\footnote{
    In general, a subspace $U$ is called $T$-invariant if
    $T(U\cap\mdef(T))\subset U$.
    If we speak of a finite-dimensional $T$-invariant subspace $U$, we
    additionally assume that $\dim U<\infty$ and $U\subset\mdef(T)$.}
  admits a basis
  consisting of Jordan chains, it is immediate that
  a Riesz basis of finite-dimensional 
  $T$-invariant subspaces is equivalent to a
  Riesz basis with parentheses of Jordan chains such that each Jordan chain
  lies inside some parenthesis.
  As a consequence of Lemma~\ref{lem_join_l2decomp},
  a Riesz basis of finite-dimensional invariant
  subspaces where almost all subspaces are one-dimensional is
  equivalent to a Riesz basis of eigenvectors and finitely many Jordan chains.
\end{remark}

\section{Spectral enclosures for $p$-subordinate perturbations}
\label{sec_psubspec}

The concept of $p$-subordination is in a certain sense an interpolation
between the notions of boundedness and relative boundedness.
We start with a result for relatively bounded perturbations.
\begin{lemma}\label{lem_generalpert}
  Let\, $G$ and $S$ be operators on a Banach space with 
  $\mdef(G)\subset\mdef(S)$ and $T=G+S$.
  If\, $0<\eps<1$ and $z\in\varrho(G)$ such that
  \begin{equation}\label{lem_generalpert-c}
    \|S(G-z)^{-1}\|\leq\eps,
  \end{equation}
  then $z\in\varrho(T)$ and
  \[\|(T-z)^{-1}\|\leq\frac{1}{1-\eps}\|(G-z)^{-1}\|
    ,\quad \|S(T-z)^{-1}\|\leq \frac{\eps}{1-\eps}.\]
  Moreover if\, $\Gamma\subset\varrho(G)$ is a simply closed,
  positively oriented integration contour and 
  (\ref{lem_generalpert-c}) holds for all
  $z\in\Gamma$, then $\Gamma\subset\varrho(T)$ and
  for the Riesz projections $Q$ and $P$ of\, $G$ and $T$ associated with 
  $\Gamma$ there are isomorphisms
  \[\range(Q)\cong\range(P),\quad \ker Q\cong\ker P.\]
\end{lemma}
\textit{Proof.}
(\ref{lem_generalpert-c}) implies the convergence of the Neumann series
\[\bigl(I+S(G-z)^{-1}\bigr)^{-1}=\sum_{k=0}^\infty\bigl(-S(G-z)^{-1}\bigr)^k\]
with
\[\big\|\bigl(I+S(G-z)^{-1}\bigr)^{-1}\big\|\leq\frac{1}{1-\|S(G-z)^{-1}\|}\leq
  \frac{1}{1-\eps}.\]
Since 
\[T-z=\bigl(I+S(G-z)^{-1}\bigr)(G-z),\]
we conclude that $z\in\varrho(T)$
with
\[\|(T-z)^{-1}\|\leq\|(G-z)^{-1}\|\big\|\bigl(I+S(G-z)^{-1}\bigr)^{-1}\big\|
  \leq\frac{1}{1-\eps}\|(G-z)^{-1}\|.\]
The identity 
$S(T-z)^{-1}=S(G-z)^{-1}(I+S(G-z)^{-1})^{-1}$
yields
$\|S(T-z)^{-1}\|\leq \eps(1-\eps)^{-1}$.

To prove the assertion about the Riesz projections, consider the operators
$T_r=G+rS$ for $r\in[0,1]$. We have the power series expansion
\[\bigl(I+rS(G-z)^{-1}\bigr)^{-1}
  =\sum_{k=0}^\infty r^k\bigl(-S(G-z)^{-1}\bigr)^k,\quad r\in[0,1],\]
which converges uniformly in $z\in\Gamma$.
Consequently $\Gamma\subset\varrho(T_r)$, and
\[(T_r-z)^{-1}=(G-z)^{-1}\bigl(I+rS(G-z)^{-1}\bigr)^{-1}\]
is continuous in $r$ uniformly for $z\in\Gamma$.
Hence the Riesz projections $P_r$ of $T_r$ associated with $\Gamma$ also
depend continuously on $r$. Now if $\|P_r-P_s\|<1$, then there are
isomorphisms
\[\range(P_r)\cong\range(P_s),\quad \ker P_r\cong\ker P_s,\]
see \cite[\S I.4.6]{kato}.
Since $r$ ranges over a compact interval, the proof is complete.
\theoremendskip

The concept of $p$-subordinate perturbations was studied by
Krein \cite[\S I.7.1]{krein} and Markus \cite[\S 5]{markus},
see also \cite[\S3.2]{wyss}.
\begin{definition}
  Let $G$, $S$ be operators on some Banach space and
  $p\in[0,1]$.
  Then $S$ is said to be \emph{$p$-subordinate to $G$}
  if $\mdef(G)\subset\mdef(S)$ and there exists $b\geq 0$ such that
  \begin{equation}\label{def_psub-b}
    \|Su\|\leq b\|u\|^{1-p}\|Gu\|^p \quad\text{for all}\quad u\in\mdef(G).
  \end{equation}
  In this case there is a minimal constant $b\geq 0$ such that 
  (\ref{def_psub-b}) holds, which is
  called the \emph{$p$-subordination bound} of $S$ to $G$.
  \definend
\end{definition}
If $S$ is $p$-subordinate to $G$ with $p<1$, then $S$ is relatively
bounded to $G$ with relative bound $0$; 
if also $0\in\varrho(G)$ and $q>p$, then
$S$ is $q$-subordinate to $G$.

\begin{remark}
  In the case that $G$ and $S$ are operators on a Hilbert space and that
  $G$ is normal with compact resolvent and $0\in\varrho(G)$, the following
  can be shown \cite[\S 5]{markus}:
  If $SG^{-p}$ is bounded with $0\leq p\leq 1$, then $S$ is $p$-subordinate
  to $G$.
  If $S$ is $p$-subordinate to $G$ with $0\leq p<1$, then
  $SG^{-q}$ is bounded for all $q>p$; in particular, $S$ is 
  relatively compact to $G$.
  \definend
\end{remark}

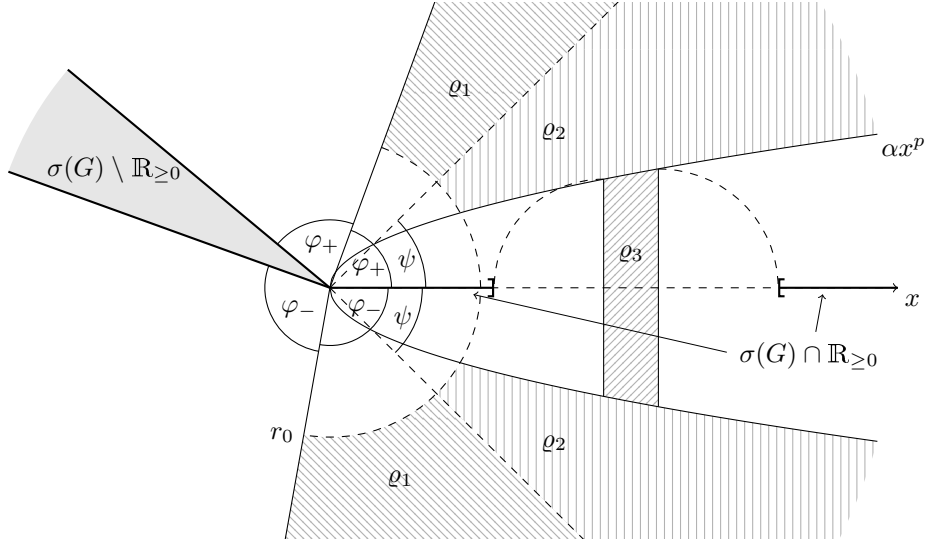
\begin{figure}
  \begin{center}
    \begin{tikzpicture}[scale=.9]
      \clip (-5,-3.7) rectangle (9,4.2);
      \ifx\tikzdontfill\undefined
      \fill[fill=black!10] (0,0) -- (140:5) arc(140:160:5) -- cycle;
      \draw (151:3.6) node {$\sigma(G)\setminus \mathbbm{R}_{\geq 0}$};
      \fill[pattern color=black!30,pattern=north east lines]
        (4,1.6)--(4,-1.6)--(4.8,-1.75)--(4.8,1.75)--cycle;
      \draw (4.4,.5) node {$\varrho_3$};
      \fill[pattern color=black!30,pattern=north west lines] 
        (45:8)--(70:6)--(70:2.2) arc(70:45:2.2)--cycle;
      \fill[pattern color=black!30,pattern=vertical lines] 
        plot[id=psubang_uf,samples=10,domain=2:8] function{.8*sqrt(x)}
        arc(16:45:8.3)--(45:8.3)--(45:2.2) arc(45:30:2.2)--cycle;
      \draw (35:4) node{$\varrho_2$} (57:3.5) node{$\varrho_1$};
      \fill[pattern color=black!30,pattern=north west lines] 
        (-45:8)--(-100:5)--(-100:2.2) arc(-100:-45:2.2)--cycle;
      \fill[pattern color=black!30,pattern=vertical lines]
        plot[id=psubang_df,samples=10,domain=2:8] function{-.8*sqrt(x)}
        arc(-16:-45:8.3)--(-45:8.3)--(-45:2.2) arc(-45:-30:2.2) --cycle;
      \draw (-35:4) node{$\varrho_2$} (-70:3) node{$\varrho_1$};

      \draw[dashed] (0,0)--(45:8);
      \draw[dashed] (0,0)--(-45:8);
      \draw[dashed,->] (0,0)--(8.3,0) node[below right=-1pt] {$x$};
      \draw (0,0)--(70:6);
      \draw (0,0)--(-100:5);
      \draw[thick] (0,0)--(140:5);
      \draw[thick] (0,0)--(-200:5);
      
      \draw[thick,[-] (2.4,0)--(0,0);
      \draw[thick,[-] (6.55,0)--(8.3,0);
      \draw[very thin,->] (7.1,-.65)--(7.2,-.1);
      \draw[very thin,->] (5.8,-.95)--(2.1,-.1);
      \draw (7,-1) node{$\sigma(G)\cap\mathbbm{R}_{\geq0}$};
	  
      \draw (4,-1.6) -- (4,1.6);
      \draw[dashed] (4,1.6) arc(90:180:1.6);
      \draw (4.8,-1.75) -- (4.8,1.75);
      \draw[dashed] (4.8,1.75) arc(90:0:1.75);

      \draw plot[id=psubang_u,samples=300,domain=0:8]
        function{.8*sqrt(x)} node[below right=-2pt] {$\alpha x^p$};
      \draw plot[id=psubang_d,samples=300,domain=0:8] function{-.8*sqrt(x)};

      \draw[dashed] (70:2.2) arc(70:-100:2.2) node[left]{$r_0$};
      \draw (.9,0) arc(0:70:.9);
      \draw (27:.65) node {$\varphi_+$};
      \draw (70:1) arc(70:140:1);
      \draw (100:.65) node {$\varphi_+$};
      \draw (.85,0) arc(0:-100:.85);
      \draw (-30:.6) node {$\varphi_-$};
      \draw (-100:.95) arc(-100:-200:.95);
      \draw (-145:.55) node {$\varphi_-$};
      \draw (1.4,0) arc(0:45:1.4);
      \draw (20:1.2) node {$\psi$};
      \draw (1.35,0) arc(0:-45:1.35);
      \draw (-22:1.15) node {$\psi$};
      \fi
    \end{tikzpicture}
  \end{center}
  \caption{The situation of Lemma~\ref{lem_psub_angle_resolvent}}
  \label{figure_psub_angle}
\end{figure}

Now we investigate how the spectrum of a normal operator $G$
changes under a $p$-subordinate perturbation $S$ with $p<1$.
We consider the case that $\sigma(G)$ lies on rays from the origin
and denote sectors in the complex plane by
\[\label{ndef_omega}
  \Omega(\varphi_-,\varphi_+)=\{re^{i\varphi}\,|\,r\geq 0\,,\,
  \varphi_-<\varphi<\varphi_+\} \quad\text{and}\quad
  \Omega(\varphi)=\Omega(-\varphi,\varphi).\]
In the next lemma the strip $\varrho_3$ corresponds to large gaps
of $\sigma(G)$ on the positive real axis, compare
Figure~\ref{figure_psub_angle}.
Sufficient conditions for the existence of such gaps may be found in
Proposition~\ref{prop_specgap}, Theorem~\ref{theo_psubpert_rootl2decomp}
and Lemma~\ref{lem_asymp2gap}.

\begin{lemma}\label{lem_psub_angle_resolvent}
  Let\, $G$ be a normal operator on a Hilbert space $H$ such that 
  $\sigma(G)\cap\Omega(2\varphi_-,2\varphi_+)\subset\bbR_{\geq 0}$
  with $-\pi\leq\varphi_-<0<\varphi_+\leq\pi$.
  Let\, $S$ be $p$-subordinate to $G$ with bound $b$,
  $0\leq p<1$, and $T=G+S$.

  Then for $\alpha>b$, $b/\alpha<\eps<1$, and 
  $0<\psi<\min\{-\varphi_-,\varphi_+,\pi/2\}$
  there exists $r_0>0$ such that the sets
  \begin{align*}
    \varrho_1&=\bigl\{z\in\overline{\Omega(\varphi_-,\varphi_+)}\,\big|\,
      |z|\geq r_0,\,z\not\in\Omega(\psi)\bigr\},\\
    \varrho_2&=\bigl\{z=x+iy\in\overline{\Omega(\psi)}\,\big|\,
      |z|\geq r_0,\,|y|\geq\alpha x^p\bigr\},\\
    \varrho_3&=\bigl\{z=x+iy\in\overline{\Omega(\psi)}\,\big|\,
      |z|\geq r_0,\,|y|\leq\alpha x^p\leq\dist(z,\sigma(G))\bigr\}
  \end{align*}
  satisfy $\varrho_1\cup\varrho_2\cup\varrho_3\subset\varrho(T)$, and
  for $z\in\varrho_1\cup\varrho_2\cup\varrho_3$ we have
  \[\|S(G-z)^{-1}\|\leq\eps,\quad
    \|(T-z)^{-1}\|\leq\frac{(1-\eps)^{-1}}{\dist(z,\sigma(G))},\quad
    \|S(T-z)^{-1}\|\leq\frac{\eps}{1-\eps}.\]
  Furthermore there is a constant\, $M>0$ such that
  \[\|(T-z)^{-1}\|\leq M\quad\text{for all}\quad
    z\in\varrho_1\cup\varrho_2\cup\varrho_3.\]
\end{lemma}
\textit{Proof.} 
We write $d=\dist(z,\sigma(G))$ and use a consequence of
the spectral theorem for normal operators,
see \cite[\S V.3.8]{kato}:
\[\|(G-z)^{-1}\|=\sup_{\lambda\in\sigma(G)}\frac{1}{|\lambda-z|}
  =\frac{1}{d},\quad
  \|G(G-z)^{-1}\|=\|I+z(G-z)^{-1}\|\leq 1+ \frac{|z|}{d}.\]
With the definition of $p$-subordination this yields
\[\|S(G-z)^{-1}u\|\leq b\|G(G-z)^{-1}u\|^p\|(G-z)^{-1}u\|^{1-p}
  \leq b\Bigl(1+\frac{|z|}{d}\Bigr)^p\frac{1}{d^{1-p}}\|u\|\]
for every $u\in H$.
In order to apply Lemma~\ref{lem_generalpert}, we thus have to show that
\begin{equation}\label{lem_psub_angle_resolvent-c}
  C=b\Bigl(1+\frac{|z|}{d}\Bigr)^p\frac{1}{d^{1-p}}\leq\eps.
\end{equation}
First we analyse the geometry of the situation:
For $z=x+iy$ we have the implications
\begin{align}
  \varphi_-\leq\arg z\leq-\frac{\pi}{2}\;\text{ or }\;
  \frac{\pi}{2}\leq\arg z\leq\varphi_+
  &\quad\Longrightarrow\quad d\geq|z|,\label{lem_psub_angle-i1}\\
  \max\Bigl\{\varphi_-,-\frac{\pi}{2}\Bigr\}\leq\arg z\leq
  \min\Bigl\{\varphi_+,\frac{\pi}{2}\Bigr\}
  &\quad\Longrightarrow\quad d\geq|y|,\label{lem_psub_angle-i2}
\end{align}
as well as
\begin{align}
  \psi\leq|\arg z|\leq\frac{\pi}{2}
  &\quad\Longrightarrow\quad |y|\geq|z|\sin\psi,
  \label{lem_psub_angle-i3}\\
  |\arg z|\leq\psi &\quad\Longrightarrow\quad x\geq|z|\cos\psi.
  \label{lem_psub_angle-i4}
\end{align}

Now let $z\in\varrho_1$. 
If $\varphi_-\leq\arg z\leq-\pi/2$ or $\pi/2\leq\arg z\leq\varphi_+$, 
then (\ref{lem_psub_angle-i1}) yields $C\leq 2^pb|z|^{p-1}\leq\eps$, 
provided $r_0$ is large enough. 
If $\psi\leq|\arg z|\leq\pi/2$, then (\ref{lem_psub_angle-i2})
and (\ref{lem_psub_angle-i3}) imply $d\geq|z|\sin\psi$ and hence
\[C\leq b\Bigl(1+\frac{1}{\sin\psi}\Bigr)^p
  \frac{1}{(|z|\sin\psi)^{1-p}}\leq\eps\]
for $r_0$ sufficiently large.

For $z\in\varrho_2$, the implications (\ref{lem_psub_angle-i2})
and (\ref{lem_psub_angle-i4}) apply and with $|y|\geq\alpha x^p$ we 
find $d\geq\alpha x^p$.
For $p>0$ we use the Minkowski inequality to get the estimate
\[\Bigl(1+\frac{|z|}{d}\Bigr)^p\leq\Bigl(1+\frac{x+|y|}{d}\Bigr)^p
  \leq1+\frac{x^p+|y|^p}{d^p}\leq 1+\frac{\alpha^{-1}d+d^p}{d^p}
  = 2+\frac{1}{\alpha}d^{1-p},\]
i.e.\ $C\leq 2bd^{p-1}+b/\alpha$.
Since $b/\alpha<\eps$ and $d\geq\alpha(|z|\cos\psi)^p$, we obtain
$C\leq\eps$ for $r_0$ sufficiently large.
On the other hand, if $p=0$ then $d\geq\alpha$ and
$C=b/d\leq b/\alpha<\eps$.

In the case $z\in\varrho_3$, (\ref{lem_psub_angle-i2}) and
(\ref{lem_psub_angle-i4}) apply, and we have $d\geq\alpha x^p$
by definition of the set $\varrho_3$.
In the same manner as for $z\in\varrho_2$,
we conclude that $C\leq\eps$ if $r_0$ is large enough.

Finally, to prove that $\|(T-z)^{-1}\|$ is uniformly bounded, we need to
show that $d^{-1}$ is bounded independently of $z$. 
For $z\in\varrho_1$ we have 
\[\text{either}\quad
  d\geq|z|\geq r_0>0 \quad\text{or}\quad 
  d\geq|z|\sin\psi\geq r_0\sin\psi>0.\]
For $z\in\varrho_2\cup\varrho_3$ we obtain
\[d\geq\alpha(|z|\cos\psi)^p\geq\alpha(r_0\cos\psi)^p>0.\]

\theoremendupskip

\begin{theo}\label{theo_psubpert_specshape}
  Let\, $G$ be a normal operator whose spectrum lies
  on finitely
  many rays $e^{i\theta_j}\bbR_{\geq 0}$ with
  $0\leq\theta_j<2\pi$, $j=1,\ldots,n$.
  Let\, $T=G+S$ where $S$ is $p$-subordinate to $G$ with bound $b$
  and $0\leq p<1$.
  Then for every $\alpha>b$ there exists $r_0>0$ such that
  \begin{equation}\label{theo_psubpert_specshape-s}
    \sigma(T)\subset B_{r_0}(0)\cup\bigcup_{j=1}^n\bigl\{
    e^{i\theta_j}(x+iy)\,\big|\,x\geq 0,\,|y|\leq\alpha x^p\bigr\},
  \end{equation}
  cf.\ Figure~\ref{figure_specshape}.
  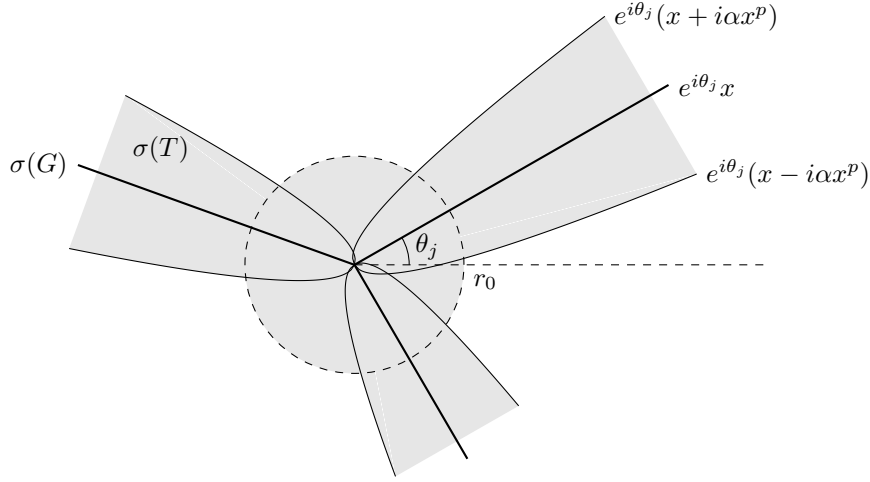
\begin{figure}
    \begin{center}
      \begin{tikzpicture}[scale=.9]
        \ifx\tikzdontfill\undefined
	\fill[fill=black!10] (1.6,0) arc(0:360:1.6);
	\fill[rotate=-60,fill=black!10]
	  plot[id=specshape_uf1,samples=100,domain=0:3] function{.6*sqrt(x)}
	  --(3,1.04)--(3,-1.04)
	  plot[id=specshape_df1,samples=100,domain=0:3] function{-.6*sqrt(x)}
	  --cycle;
	
	\fill[rotate=160,fill=black!10]
	  plot[id=specshape_uf2,samples=100,domain=0:4] function{.6*sqrt(x)}
	  --(4,1.2)--(4,-1.2)
	  plot[id=specshape_df2,samples=100,domain=0:4] function{-.6*sqrt(x)}
	  --cycle;
	
	\fill[rotate=30,fill=black!10]
	  plot[id=specshape_uf3,samples=100,domain=0:5] function{.6*sqrt(x)}
	  --(5,1.34)--(5,-1.34)
	  plot[id=specshape_df3,samples=100,domain=0:5] function{-.6*sqrt(x)}
	  --cycle;
	\fi

	\draw[dashed] (0,0)--(6,0);
	
	\draw[thick] (0,0)--(-60:3.3);
	\draw[rotate=-60] plot[id=specshape_u1,samples=200,domain=0:3]
	  function{.6*sqrt(x)};
	\draw[rotate=-60] plot[id=specshape_d1,samples=200,domain=0:3] 
	  function{-.6*sqrt(x)};
	
	\draw[thick] (0,0)--(160:4.3) node[left]{$\sigma(G)$};
        \ifx\tikzdontfill\undefined
	\draw[rotate=160] plot[id=specshape_u2,samples=200,domain=0:4]
	  function{.6*sqrt(x)};
	\draw[rotate=160] plot[id=specshape_d2,samples=200,domain=0:4] 
	  function{-.6*sqrt(x)};
	\draw (149:3.3) node{$\sigma(T)$};
	
	\draw[thick] (0,0)--(30:5.3) node[right]{$e^{i\theta_j}x$};
	\fi
	\draw[rotate=30] plot[id=specshape_u3,samples=200,domain=0:5]
	  function{.6*sqrt(x)} node[right]{$e^{i\theta_j}(x+i\alpha x^p)$};
	\draw[rotate=30] plot[id=specshape_d3,samples=200,domain=0:5] 
	  function{-.6*sqrt(x)} node[right]{$e^{i\theta_j}(x-i\alpha x^p)$};
	
        \ifx\tikzdontfill\undefined
	\draw[dashed] (1.6,0) node[below right]{$r_0$} arc(0:360:1.6);

	\draw (.8,0) arc(0:30:.8);
	\draw (14:1.1) node{$\theta_j$};
	\fi
      \end{tikzpicture}
    \end{center}
    \caption{The spectrum after a $p$-subordinate perturbation}
    \label{figure_specshape}
  \end{figure}
  If\, $G$ has a compact resolvent, then so has $T$.
\end{theo}
\textit{Proof.}
Without loss of generality, we assume 
$\theta_1<\theta_2<\ldots<\theta_n$ and set $\theta_0=\theta_n-2\pi$,
$\theta_{n+1}=\theta_0+2\pi$. Then we may, after a rotation by $\theta_j$,
apply Lemma~\ref{lem_psub_angle_resolvent} to each sector
$\Omega(\theta_{j-1},\theta_{j+1})$. More precisely, we apply the lemma
to the operators $e^{-i\theta_j}G$, $e^{-i\theta_j}S$, $e^{-i\theta_j}T$
with $\varphi_+=(\theta_{j+1}-\theta_j)/2$, 
$\varphi_-=(\theta_{j-1}-\theta_j)/2$, and some suitable $\eps$. 
This yields the implication
\begin{align*}
  &z\in\sigma(T),\quad\frac{\theta_{j-1}+\theta_j}{2}\leq\arg z\leq
  \frac{\theta_j+\theta_{j+1}}{2},\quad|z|\geq r_0 \\
  &\Longrightarrow\quad
  z\in\{e^{i\theta_j}(x+iy)\,\big|\,x\geq 0,\,|y|\leq\alpha x^p\bigr\}
\end{align*}
with some $r_0\geq 0$ for each $j=1,\ldots,n$.
If $G$ has compact resolvent, the identity
\[(T-z)^{-1}=(G-z)^{-1}\bigl(I+S(G-z)^{-1}\bigr)^{-1}
  \quad\text{for}\quad z\in\varrho(G)\cap\varrho(T)\]
implies that $T$ has compact resolvent too.
\theoremendskip

The statement about the asymptotic shape of the spectrum of 
$T$ can be refined as follows:
\begin{remark}\label{rem_psubpert_specshape}
  To obtain a condition for $z\in\varrho(T)$, we consider without loss
  of generality the case $\sigma(G)\cap\Omega(2\varphi)\subset\bbR_{\geq 0}$,
  $0<\varphi\leq\pi/2$, and $z=x+iy\in\overline{\Omega(\varphi)}$.
  Then $\dist(z,\sigma(G))\geq|y|$ and,
  in view of (\ref{lem_psub_angle_resolvent-c}),
  $b(1+|z|/|y|)^p|y|^{p-1}<1$ is sufficient to get $z\in\varrho(T)$.
  For $p>0$ this leads to the condition
  \[x<\left(\frac{|y|}{b}\right)^{1/p}
    \sqrt{1-2b^{1/p}|y|^{1-1/p}},\]
  which is asymptotically better than $x<(|y|/\alpha)^{1/p}$ from the theorem
  since $1-2b^{1/p}|y|^{1-1/p}\to 1$ as $|y|\to\infty$.
  For $p=0$ we obtain the optimal condition $b<|y|$.

  For $p>0$, the estimates of Markus~\cite[Lemma~5.2]{markus}
  lead to asymptotics which are even slightly better.
  Also note that simply taking the limit $\alpha\to b$ in 
  Theorem~\ref{theo_psubpert_specshape}
  is not possible since then also $r_0\to\infty$.
  \definend
\end{remark}

\section{Estimates for Riesz projections}
\label{sec_rprojest}

In this section $G$ is always a normal operator with compact
resolvent on a Hilbert space $H$ such that
\[\sigma(G)\cap\Omega(2\varphi)\subset\bbR_{\geq0}
  \quad\text{with}\quad 0<\varphi\leq\frac{\pi}{2}\]
and $T=G+S$ with $S$ $p$-subordinate to $G$ and $p<1$.

The first two lemmas can be found in the book of Markus \cite{markus},
for the special case $\alpha=4b$.
Since his proofs literally apply to the general situation, we omit them here;
see also \cite[\S3.3]{wyss}.
\begin{lemma}[\mbox{Markus \cite[Lemma 6.6]{markus}}]
  \label{lem_estim_contour1}
  Let\, $G$ be normal with compact resolvent and
  $\sigma(G)\cap\Omega(2\varphi)\subset\bbR_{\geq 0}$
  with $0<\varphi\leq\pi/2$.
  Then for $0\leq p<1$, $\alpha>0$
  there exists $r_0>0$ such that
  the contours
  \begin{equation}\label{lem_estim_contour1-contdef}
    \Gamma_\pm=\{x+iy\in\bbC\,|\,x\geq r_0,\,y=\pm \alpha x^p\}
  \end{equation}
  satisfy $\Gamma_\pm\subset\varrho(G)\cap\overline{\Omega(\varphi)}$
  and we have
  \[\int_{\Gamma_\pm}|z|^p\|(G-z)^{-1}u\|^2\,|dz|\leq C_1\|u\|^2,\quad
    \int_{\Gamma_\pm}|z|^{p-2}\|G(G-z)^{-1}u\|^2\,|dz|\leq C_2\|u\|^2\]
  for all\, $u\in H$ with some constants $C_1,C_2\geq 0$.
  \theoremend
\end{lemma}

\begin{lemma}[\mbox{Markus \cite[Lemma 6.7]{markus}}]
  \label{lem_estim_contour2}
  Let\, $G$ be normal with compact resolvent and
  $\sigma(G)\cap\Omega(2\varphi)\subset\bbR_{\geq 0}$
  with $0<\varphi\leq\pi/2$.
  Let\, $(x_k)_{k\geq1}$ be a sequence of positive numbers,
  $0\leq p<1$, and $\alpha,c_1,c_2>0$ such that\,
  $\alpha x_1^{p-1}\leq\tan\varphi$ and
  \[x_n^{1-p}-x_k^{1-p}\geq c_1(n-k)\quad\text{for}\quad n>k,\quad
    \dist(x_k,\sigma(G))\geq c_2x_k^p \quad\text{for}\quad k\geq 1.\]
  Then the lines
  \begin{equation}\label{lem_estim_contour2-contdef}
    \gamma_k=\bigl\{x_k+iy\in\bbC\,\big|\,|y|\leq\alpha x_k^p\bigr\}
  \end{equation}
  satisfy $\gamma_k\subset\varrho(G)\cap\overline{\Omega(\varphi)}$ and we have
  \[\sum_{k=1}^\infty x_k^p\int_{\gamma_k}\|(G-z)^{-1}u\|^2\,|dz|
    \leq C_1\|u\|^2,\quad
    \sum_{k=1}^\infty x_k^{p-2}\int_{\gamma_k}\|G(G-z)^{-1}u\|^2\,|dz|
    \leq C_2\|u\|^2\]
  for all\, $u\in H$ with some constants $C_1,C_2\geq 0$.
  \theoremend
\end{lemma}

With the previous resolvent estimates at hand, we derive an estimate for
a sequence of Riesz projections associated with the parabola $\Gamma_\pm$
and the lines $\gamma_k$:
\begin{lemma}\label{lem_rieszproj_estim}
  Let\, $G$ be normal with compact resolvent,
  $\sigma(G)\cap\Omega(2\varphi)\subset\bbR_{\geq 0}$ 
  with $0<\varphi\leq\pi/2$,
  $S$ $p$-subordinate to $G$ with bound $b$,
  $0\leq p<1$, and $T=G+S$.
 
  Let\, $\alpha>b$, let $(x_k)_{k\geq 1}$, $\gamma_k$ be as in
  Lemma~\ref{lem_estim_contour2}, and suppose that there is a constant
  $M\geq 0$ such that
  \[\gamma_k\subset\varrho(T) \quad\text{and}\quad
    \|S(T-z)^{-1}\|\leq M
    \quad\text{for all}\quad z\in\gamma_k,\,k\geq 1.\]
  Then there exist\, $r_0>0$, $k_0\geq 1$ such that\, $x_{k_0}\geq r_0$ and 
  the following holds:
  If\, $\Gamma_\pm$ is as in (\ref{lem_estim_contour1-contdef})
  and $\Gamma_k$ with $k\geq k_0$ is the positively oriented boundary contour 
  of the region
  enclosed by $\gamma_k,\Gamma_-,\gamma_{k+1},\Gamma_+$,
  then $\Gamma_k\subset\varrho(T)$.
  If\, $P_k$ is the Riesz projection of\, $T$
  associated with $\Gamma_k$, then
  \begin{equation}\label{lem_rieszproj_estim-e}
    \sum_{k=k_0}^\infty|(P_ku|v)|\leq C\|u\|\|v\|\quad\text{for all}\quad
    u,v\in H
  \end{equation}
  with some constant\, $C\geq 0$.
\end{lemma}
\textit{Proof.}
We want to apply Lemmas~\ref{lem_psub_angle_resolvent}, 
\ref{lem_estim_contour1} and \ref{lem_estim_contour2},
and choose $\eps\in\,]b/\alpha,1[$ and $r_0$ accordingly. The assumptions
on the sequence $(x_k)_k$ imply that it tends monotonically to infinity
and we choose $k_0$ such that $x_{k_0}\geq r_0$.
By Lemma~\ref{lem_psub_angle_resolvent}, $\|S(T-z)^{-1}\|$ is uniformly
bounded on $\Gamma_\pm$. We thus have
\[\Gamma_k\subset\varrho(G)\cap\varrho(T)\quad\text{and}\quad
  \|S(T-z)^{-1}\|\leq M_0\quad\text{for all}\quad z\in\Gamma_k,k\geq k_0,\]
with some $M_0\geq 0$.
Consider now the Riesz projections $Q_k$ of $G$ associated with $\Gamma_k$,
which are orthogonal since $G$ is normal. 
It is easy to see that, to prove (\ref{lem_rieszproj_estim-e}),
it suffices to prove
\[\sum_{k=k_0}^\infty\big|\big((P_k-Q_k)u\big|v\big)\big|\leq C\|u\|\|v\|.\]
Now
\[P_k-Q_k=\frac{i}{2\pi}\int_{\Gamma_k}\left((T-z)^{-1}-(G-z)^{-1}\right)dz
  =\frac{-i}{2\pi}\int_{\Gamma_k}(T-z)^{-1}S(G-z)^{-1}dz\]
and hence
\[\big|\big((P_k-Q_k)u\big|v\big)\big|\leq\frac{1}{2\pi}\int_{\Gamma_k}
  \|S(G-z)^{-1}u\|\|(T-z)^{-*}v\|\,|dz|.\]
Then, with the help of
\begin{align*}
  (T-z)^{-1}&=(G-z)^{-1}\left(I-S(T-z)^{-1}\right)\\
  \Longrightarrow\qquad \|(T-z)^{-*}v\|&\leq\big(1+
  \underbrace{\|S(T-z)^{-1}\|}_{\leq M_0}\big)\|(G-z)^{-*}v\|
\end{align*}
and $\|(G-z)^{-*}v\|=\|(G-z)^{-1}v\|$ (since $G$ is normal), we find
\begin{align*}
  \sum_{k=k_0}^\infty\big|\big((&P_k-Q_k)u\big|v\big)\big|
  \leq\frac{1+M_0}{2\pi}\sum_{k=k_0}^\infty\int_{\Gamma_k}
  \|S(G-z)^{-1}u\|\|(G-z)^{-1}v\|\,|dz|\\
  &\leq\frac{1+M_0}{2\pi}\left(\int_{\Gamma_+}+\int_{\Gamma_-}
  +2\sum_{k=k_0}^\infty\int_{\gamma_k}\right)
    \|S(G-z)^{-1}u\|\|(G-z)^{-1}v\|\,|dz|.
\end{align*}
Using $p$-subordination, Lemma~\ref{lem_estim_contour1},
and (for $p\neq0$) H\"older's inequality, we estimate
\[\begin{split}
    \int_{\Gamma_\pm}\|&S(G-z)^{-1}u\|\|(G-z)^{-1}v\|\,|dz|\\ 
    &\leq\biggl(
    \int_{\Gamma_\pm}|z|^{-p}\|S(G-z)^{-1}u\|^2\,|dz|\biggr)^{1/2}
    \biggl(\underbrace{\int_{\Gamma_\pm}|z|^p\|(G-z)^{-1}v\|^2\,|dz|}
      _{\leq C_1\|v\|^2} \biggr)^{1/2},
  \end{split} \]
\begin{align*}
  \int_{\Gamma_\pm}&|z|^{-p}\|S(G-z)^{-1}u\|^2\,|dz|\\
  &\leq b^2\left(\int_{\Gamma_\pm}|z|^{p-2}\|G(G-z)^{-1}u\|^2\,|dz|\right)^p
  \left(\int_{\Gamma_\pm}|z|^p\|(G-z)^{-1}u\|^2\,|dz|\right)^{1-p}\\
  &\leq b^2C_2^pC_1^{1-p}\|u\|^2.
\end{align*}
In the same way, with Lemma~\ref{lem_estim_contour2}, we see that
\begin{align*}
  \sum_k&\int_{\gamma_k}\|S(G-z)^{-1}u\|\|(G-z)^{-1}v\|\,|dz|\\
  &\leq\left(\sum_k\int_{\gamma_k}x_k^{-p}\|S(G-z)^{-1}u\|^2\,|dz|
  \right)^{1/2}
  \Biggl(\underbrace{\sum_k\int_{\gamma_k}x_k^p
   \|(G-z)^{-1}v\|^2\,|dz|}_{\leq C_1^\prime\|v\|^2}\Biggr)^{1/2}
\end{align*}
and
\begin{align*}
  \sum_k&\int_{\gamma_k}x_k^{-p}\|S(G-z)^{-1}u\|^2\,|dz|\\
  &\leq b^2\left(\sum_k\int_{\gamma_k}x_k^{p-2}\|G(G-z)^{-1}u\|^2
  \,|dz|\right)^p
  \left(\sum_k\int_{\gamma_k}x_k^p\|(G-z)^{-1}u\|^2\,|dz|\right)^{1-p}\\
  &\leq b^2C_2^{\prime p}C_1^{\prime1-p}\|u\|^2.
\end{align*}

\theoremendupskip

To proceed, we need the concept of 
the determinant for operators, see 
\cite[\S 2.5]{markus},
\cite[Chapter VII]{gohberg-goldberg-kaashoek} and
\cite[\S IV.1]{gohberg-krein}.
For an operator $A$ of finite rank $m$,
the \emph{determinant} of $I+A$ is defined by
\begin{equation}
  \det(I+A)=\det\bigl((I+A)|_{\range(A)}\bigr)
\end{equation}
and it satisfies
\begin{itemize}
\item[(i)] $|\det(I+A)|\leq(1+\|A\|)^m$;
\item[(ii)] $I+A$ is invertible if and only if\, $\det(I+A)\neq0$, and
  in this case
  \[\|(I+A)^{-1}\|\leq\frac{(1+\|A\|)^m}{|\det(I+A)|}\,;\]
\item[(iii)] if the operator-valued function $B:\Omega\to L(H)$ is 
  analytic on a domain $\Omega\subset\bbC$, then 
  \(z\mapsto\det(I+AB(z))\)
  is analytic on $\Omega$ too.
\end{itemize}

We also use the following  auxiliary result from complex analysis,
cf.\ \cite[Lemma~1.6]{markus}, \cite[Theorem~I.11]{levin}:
\begin{lemma}\label{lem_func_discs}
  Let\, $U\subset\bbC$ be a bounded, simply connected domain, $F\subset U$
  compact, $z_0$ an interior point of\, $F$, and $\eta>0$.
  Then there exists a constant\, $C>0$
  such that the following holds: If\, $a,b\in\bbC$ and $f:aU+b\to\bbC$ 
  with $f(az_0+b)\neq0$ is 
  holomorphic and bounded, then there is a set\, $E\subset\bbC$ being the union
  of finitely many discs with radii summing up to at most\, $|a|\eta$
  such that
  \[|f(z)|\geq\frac{|f(az_0+b)|^{1+C}}{\|f\|_{aU+b,\infty}^C}
    \quad\text{for all}\quad z\in(aF+b)\setminus E.\]
\end{lemma}

The next proposition permits us to estimate the resolvent
of the perturbed operator even close to its eigenvalues
by artificially creating a gap in the spectrum of $G$.
The method is taken from Lemma~5.6 in \cite{markus},
which may be obtained from our proposition as a corollary.
We denote by $N_+(r_1,r_2,G)$ the sum of the multiplicities of all 
the eigenvalues of $G$ in the open interval $]r_1,r_2[$,
\begin{equation}\label{eigmultsum_int}
  N_+(r_1,r_2,G)=\!\!\sum_{\lambda\in\sigma_p(G)\cap\,]r_1,r_2[}\!\!
  \dim \rsub(\lambda).
\end{equation}
\begin{prop}\label{prop_resolvent_discs}
  Let\, $G$ be normal with compact resolvent,
  $\sigma(G)\cap\Omega(2\varphi)\subset\bbR_{\geq0}$ with
  $0<\varphi\leq\pi/2$, 
  $S$ $p$-subordinate to $G$ with bound $b$, $0\leq p<1$,
  and $T=G+S$.

  Let\, $l>b$, $0\leq l_0<l-b$ and $\eta>0$. Then there are constants
  $C_0,C_1,r_0>0$ such that for every $r\geq r_0$ there is a set
  $E_r\subset\bbC$ with the following properties:
  \begin{itemize}
  \item[(i)] $E_r$ is the union of finitely many discs with radii summing 
    up to at most\, $\eta r^p$.
  \item[(ii)] For every $z\in\overline{\Omega(\varphi)}\setminus E_r$ 
    with $|\Real z-r|\leq l_0r^p$ we have 
    \[z\in\varrho(T)\quad\text{and}\quad
      \|(T-z)^{-1}\|\leq\frac{C_0C_1^m}{r^p},\quad
      \|S(T-z)^{-1}\|\leq C_0C_1^m\]
    where $m=N_+(r-lr^p,r+lr^p,G)$.
  \end{itemize}
\end{prop}
\textit{Proof.}
We choose $l_1\in\,]l_0,l-b[$ and $\alpha$, $\tilde{b}$ such that
\[b<\tilde{b}<\alpha<l-l_1.\]
Let $r\geq r_0$. We may assume that $r-lr^p>0$ by choosing $r_0$ large enough.
Let $\lambda_1,\ldots,\lambda_n$ be the eigenvalues of $G$ in 
$\Delta_r=\,]r-lr^p,r+lr^p[\,$,
$P_1,\ldots,P_n$ the orthogonal projections onto the corresponding 
eigenspaces, and
\[K_r=\sum_{j=1}^n(\lambda_j-\widetilde{\lambda}_j)P_j
  \quad\text{with}\quad
  \widetilde{\lambda}_j=\begin{cases}r-lr^p&\text{if }\,\lambda_j<r,\\
  r+lr^p&\text{if }\,\lambda_j\geq r.\end{cases}\]
Then $G_r=G-K_r$ is a normal operator with 
$\sigma(G_r)\cap\Omega(2\varphi)\subset\bbR_{\geq0}$ and 
$\Delta_r\subset\varrho(G_r)$. $K_r$ has rank $m$ and satisfies
$\|K_r\|\leq lr^p$. 
Noting that
$\lambda_j/\widetilde{\lambda}_j\leq r/(r-lr^p)$ for all $j$,
it is straightforward to show that
\[\|Gu\|\leq\frac{r}{r-lr^p}\|G_ru\|.\]
Since $1-lr^{p-1}\to 1$ as $r\to\infty$ and $b<\tilde{b}$,
we conclude
\[\|Su\|\leq b\|Gu\|^p\|u\|^{1-p}\leq b\Bigl(\frac{1}{1-lr^{p-1}}\Bigr)^p
  \|G_ru\|^p\|u\|^{1-p}\leq\tilde{b}\|G_ru\|^p\|u\|^{1-p},\]
provided $r_0$ is sufficiently large. Thus $S$ is $p$-subordinate
to $G_r$ with bound less or equal than $\tilde{b}$.

Next, we want to prove that
\begin{equation}\label{resolvent_discs_s_incl}
  |x-r|\leq l_1r^p \quad\Rightarrow\quad 
  ]x-\alpha x^p,x+\alpha x^p[\,\subset\varrho(G_r)
\end{equation}
for $r_0$ sufficiently large.
Let $|x-r|\leq l_1r^p$. Since the function $x\mapsto x-\alpha x^p$ is
monotonically increasing for large $x$, we have 
\[x-\alpha x^p\geq r-l_1r^p-\alpha \bigl(r-l_1r^p\bigr)^p
  \geq r-l_1r^p-\alpha r^p >r-lr^p\]
for $r_0$ large enough. Furthermore
$\alpha(1+l_1r^{p-1})^p\leq l-l_1$ holds for large $r$ and we obtain
\[x+\alpha x^p\leq r+l_1r^p+\alpha\bigl(r+l_1r^p\bigr)^p \leq r+lr^p.\]
In view of $\Delta_r\subset\varrho(G_r)$, (\ref{resolvent_discs_s_incl})
is proved.

Now we aim to apply Lemma~\ref{lem_func_discs} to the function
\[d(z)=\det(I+K_r(G_r+S-z)^{-1}), \quad z\in\varrho(G_r+S),\]
and the sets
\begin{align*}
  U_r&=\bigl\{x+iy\,\big|\,|x-r|<l_1r^p,\,|y|<4br^p\bigr\},\\
  F_r&=\bigl\{x+iy\,\big|\,|x-r|\leq l_0r^p,\,
  |y|\leq 3br^p\bigr\}.
\end{align*}
For $r_0$ sufficiently large we have $U_r\subset\Omega(\varphi)$.
Using (\ref{resolvent_discs_s_incl}),
we can apply Lemma~\ref{lem_psub_angle_resolvent} to $G_r+S$ with some
$\eps\in\,]\tilde{b}/\alpha,1[$; we obtain $U_r\subset\varrho(G_r+S)$ and,
for $z\in U_r$,
\[\dist(z,\sigma(G_r))\geq lr^p-l_1r^p>\alpha r^p\]
and
\[\|(G_r+S-z)^{-1}\|
  \leq\frac{(1-\eps)^{-1}}{\alpha r^p},\;\;\|S(G_r+S-z)^{-1}\|\leq
  \frac{\eps}{1-\eps}.\]
Then
\[|d(z)|\leq\left(1+\|K_r\|\|(G_r+S-z)^{-1}\|\right)^m
  \leq\left(1+\frac{l(1-\eps)^{-1}}{\alpha}\right)^m=c_0^m\]
on $U_r$ with $c_0>0$.
For $z\in\varrho(T)\cap U_r$ we have
\[I=\left(I+K_r(G_r+S-z)^{-1}\right)\left(I-K_r(T-z)^{-1}\right).\]
Applying Lemma~\ref{lem_psub_angle_resolvent} (now with $\eps=2/3$) to 
the operator $T$ and $z_r=r+i\cdot 2br^p\in F_r$, we obtain
\[z_r\in\varrho(T)\quad\text{and}\quad
  \|(T-z_r)^{-1}\|\leq\frac{3}{2br^p}\]
and thus
\[\left|\frac{1}{d(z_r)}\right|=\left|\det\left(I-K_r(T-z_r)^{-1}\right)
  \right|\leq\left(1+\frac{3l}{2b}\right)^m=c_1^m\]
with $c_1>0$.
Lemma~\ref{lem_func_discs} then yields
a constant $C>0$
depending only on $b,l_0,l_1$ and $\eta$ such that for every $r\geq r_0$
there exists a union $E_r$ of discs with radii summing up to at most 
$\eta r^p$ and
\[|d(z)|\geq c_0^{-mC}c_1^{-m(1+C)} \quad\text{for all}\quad
  z\in F_r\setminus E_r.\]
For $z\in F_r\setminus E_r$, we thus obtain that
$I+K_r(G_r+S-z)^{-1}$ is invertible with
\[\big\|\bigl(I+K_r(G_r+S-z)^{-1}\bigr)^{-1}\big\|
  \leq\frac{c_0^m}{|d(z)|}\leq(c_0c_1)^{(1+C)m}.\]
Consequently $z\in\varrho(T)$ with
\[(T-z)^{-1}=(G_r+S-z)^{-1}\bigl(I+K_r(G_r+S-z)^{-1}\bigr)^{-1}\]
and
\begin{align*}
  \|(T-z)^{-1}\|&\leq\frac{(1-\eps)^{-1}}{\alpha r^p}(c_0c_1)^{(1+C)m}
  \leq\frac{C_0C_1^m}{r^p},\\
  \|S(T-z)^{-1}\|&\leq\frac{\eps}{1-\eps}(c_0c_1)^{(1+C)m}
  \leq C_0C_1^m
\end{align*}
with appropriate constants $C_0,C_1$ depending on 
$b,l,l_0,l_1,\eta,\alpha,\eps$ only.

Finally, we consider $z=x+iy\in\overline{\Omega(\varphi)}$ with
$|x-r|\leq l_0r^p$ and $|y|\geq 3br^p$.
Then
\[2bx^p\leq 2b\left(r+l_0r^p\right)^p
  \leq 3br^p\leq|y|\]
holds for $r_0$ sufficiently large.
Applying Lemma~\ref{lem_psub_angle_resolvent} (again with $\eps=2/3$),
we obtain $z\in\varrho(T)$ and
\[\|(T-z)^{-1}\|\leq\frac{3}{|y|}\leq\frac{1}{br^p}\leq\frac{C_0C_1^m}{r^p},
  \quad \|S(T-z)^{-1}\|\leq 2\leq C_0C_1^m\]
for $C_0\geq\max\{2,b^{-1}\}$ and $C_1\geq1$.
\theoremend

\begin{coroll}\label{coroll_resolvent_discs2}
  Let\, $G$ be normal with compact resolvent,
  $\sigma(G)\cap\Omega(2\varphi)\subset\bbR_{\geq0}$ with
  $0<\varphi\leq\pi/2$, 
  $S$ $p$-subordinate to $G$ with bound $b$, $0\leq p<1$,
  and $T=G+S$.

  Then for $l_0,q>0$ there are constants
  $C_0,C_1,r_0>0$ such that for every $r\geq r_0$ the following holds:
  For every $z=x+iy$ with $|x-r|\leq l_0r^p$, $|y|\leq 2bx^p$
  there exists $q_1\in\,]0,q[$ such that
  \[|w-z|=q_1r^p \quad\Longrightarrow\quad w\in\varrho(T),\quad 
    \|(T-w)^{-1}\|\leq \frac{C_0C_1^m}{r^{p}},\]
  where $m=N_+(r-lr^p,r+lr^p,G)$ with $l=b+2(l_0+q)$.
\end{coroll}
\textit{Proof.}
We use Proposition~\ref{prop_resolvent_discs} with 
$l=b+2(l_0+q)$, $l_0+q$ replacing $l_0$, and $\eta=q/3$.
For $z$ as in the claim and $|w-z|\leq qr^p$ we have $|\arg w|\leq\varphi$
(for $r_0$ large enough) and
\(|\Real w-r|\leq (l_0+q)r^p.\)
Now the sum of the diameters of the discs in $E_r$ is at most
$2\eta r^p<qr^p$. Hence there exists $q_1\in]0,q[$ such that
$w\not\in E_r$ for $|w-z|=q_1r^p$ and the claim is proved.
\theoremendskip

Under certain assumptions on the distribution of the eigenvalues of $G$
on the positive real axis, we now obtain a sequence of closed integration 
contours in $\varrho(T)$ of the form in Lemma~\ref{lem_rieszproj_estim}
and estimates for the associated Riesz projections.
\begin{prop}\label{prop_rieszproj_estim}
  Let\, $G$ be normal with compact resolvent,
  $\sigma(G)\cap\Omega(2\varphi)\subset\bbR_{\geq0}$ with
  $0<\varphi\leq\pi/2$, 
  $S$ $p$-subordinate to $G$ with bound $b$, $0\leq p<1$,
  and $T=G+S$.

  Assume that there is a sequence $(r_k)_{k\geq 1}$ of positive numbers
  tending monotonically to infinity and some $l>b$, $m\in\bbN_{\geq1}$
  such that
  \begin{equation}
    N_+(r_k-lr_k^p,r_k+lr_k^p,G)\leq m\quad\text{for all}\quad k\geq 1.
  \end{equation}
  Then there are constants $C,r_0>0$, $\alpha>b$, and a sequence 
  $(x_k)_{k\geq 1}$ in $\bbR_{\geq0}$
  tending monotonically to infinity such that the following holds:
  \begin{itemize}
  \item[(i)] $z\in\overline{\Omega(\varphi)}$ with $\Real z=x_k$ implies
    $z\in\varrho(T),\,\|(T-z)^{-1}\|\leq C$.
  \item[(ii)] The contours 
    $\Gamma_\pm,\gamma_k$ from (\ref{lem_estim_contour1-contdef}) and
    (\ref{lem_estim_contour2-contdef}) satisfy
    $\Gamma_\pm,\gamma_k\subset\varrho(T)$.
  \item[(iii)] If\, $P_k$ is the Riesz projection of\, $T$
    associated with the region
    enclosed by $\gamma_k,\Gamma_-,\gamma_{k+1},\Gamma_+$, then
    \[\sum_{k=1}^\infty|(P_ku|v)|\leq C\|u\|\|v\|\quad\text{for all}\quad
      u,v\in H.\]
  \end{itemize}
\end{prop}
\textit{Proof.}
We apply Proposition~\ref{prop_resolvent_discs} with
$l_0=(l-b)/2$ and $\eta=l_0/2$ to $r=r_k$, $k\geq k_0$, $k_0$
appropriate. Since the sum of the diameters of the discs in $E_r$
is at most $l_0r_k^p$ and the interval $[r_k-l_0r_k^p,r_k+l_0r_k^p]$
contains at most $m$ eigenvalues of $G$, we can find an $x_k$ such that
\[|x_k-r_k|\leq l_0r_k^p, \quad
  \dist(x_k,\sigma(G))\geq\frac{l_0}{3m}r_k^p,\]
and that $z\in\overline{\Omega(\varphi)}$ with $\Real z=x_k$ implies
\[z\in\varrho(T),\quad
  \|(T-z)^{-1}\|\leq \frac{C_0C_1^m}{r_k^{p}},\quad 
  \|S(T-z)^{-1}\|\leq C_0C_1^m.\]
Then $x_k/r_k\to 1$ as $k\to\infty$ and we obtain
\[\dist(x_k,\sigma(G))\geq c_2x_k^p \quad\text{for}\quad k\geq k_0\]
with $c_2>0$ and $k_0$ appropriately chosen.
Since $x_k\to\infty$, for every $k_1$ there exists $k_2>k_1$ such that
$x_{k_2}^{1-p}-x_{k_1}^{1-p}\geq 1$.
Passing to an appropriate subsequence, we can thus assume that
\[x_{k+1}^{1-p}-x_k^{1-p}\geq 1\quad\text{for all}\quad k\in\bbN,\]
which yields
\[x_n^{1-p}-x_k^{1-p}\geq n-k \quad\text{for}\quad n>k.\]
Now an application of Lemma~\ref{lem_rieszproj_estim} with
$\alpha=2b$ and
the sequence $(x_k)_{k\geq k_0}$, $k_0$ large enough, completes
the proof.
\theoremendskip

If the spectrum of $G$ has sufficiently large gaps on $\bbR_{\geq 0}$,
then the spectrum of $T$ has corresponding gaps 
(cf.\ Figure~\ref{figure_projdist}):
\begin{figure}
  \begin{center}
    \begin{tikzpicture}[scale=.9]
      \ifx\tikzdontfill\undefined
      \fill[pattern color=black!30,pattern=north west lines] 
        plot[id=projdist_uf3,samples=20,domain=4.5:6.9] function{.8*sqrt(x)}
        --(6.9,2.1)--(6.9,-2.1)
        --plot[id=projdist_df3,samples=20,domain=4.5:6.9] function{-.8*sqrt(x)}
        --(4.5,-1.7)--(4.5,1.7)--cycle;
      \fi

      \draw[->] (0,0)--(9.8,0) node[below right=-1pt]{$x$};
      \draw[dashed] plot[id=projdist_ud,samples=50,domain=0:1]
        function{.8*sqrt(x)};
      \draw[dashed] plot[id=projdist_dd,samples=50,domain=0:1]
        function{-.8*sqrt(x)};
      \draw plot[id=projdist_u,samples=100,domain=1:9.5]
        function{.8*sqrt(x)} node[right] {$\Gamma_+$};
      \draw plot[id=projdist_d,samples=100,domain=1:9.5] function{-.8*sqrt(x)}
        node[right]{$\Gamma_-$};

      \ifx\tikzdontfill\undefined
      \draw[dashed] (1,-.8)--(1,.8);
      \draw (1,0) node[below right]{$r_0$};
      
      \draw (2.5,-1.26)--(2.5,1.26);
      \draw (4.5,-1.7)--(4.5,1.7);
      \draw[thick,[-] (3.17,0)--(3.83,0);
      \draw[thick,[-] (3.83,0)--(3.8,0);
      \draw (6.9,-2.1)--(6.9,2.1);
      \draw (9,-2.4)--(9,2.4);
      \draw[thick,[-] (7.6,0)--(8.3,0);
      \draw[thick,[-] (8.3,0)--(8,0);

      \draw (2.5,.6) node[left]{$\gamma_k^-$};
      \draw (4.5,.85) node[left]{$\gamma_k^+$};
      \draw (6.9,1) node[right]{$\gamma_{k+1}^-$};
      \draw (9,1.2) node[right]{$\gamma_{k+1}^+$};
      \draw (5.7,.6) node{$\varrho(T)$};
      \draw (7.8,-1.2) node{$\sigma(G)$};
      \draw[very thin,->] (7.2,-1.15)--(3.9,-.15);
      \draw[very thin,->] (7.8,-.8)--(7.95,-.1);
      \fi
      \draw (3.5,-2.5) node{$[r_k-\beta r_k^p,r_k+\beta r_k^p]$};
      \draw[very thin,->] (3.5,-2.1)--(3.5,-.1);
    \end{tikzpicture}
  \end{center}
  \caption{A large gap in $\sigma(G)$ yields a gap in $\sigma(T)$.}
  \label{figure_projdist}
\end{figure}
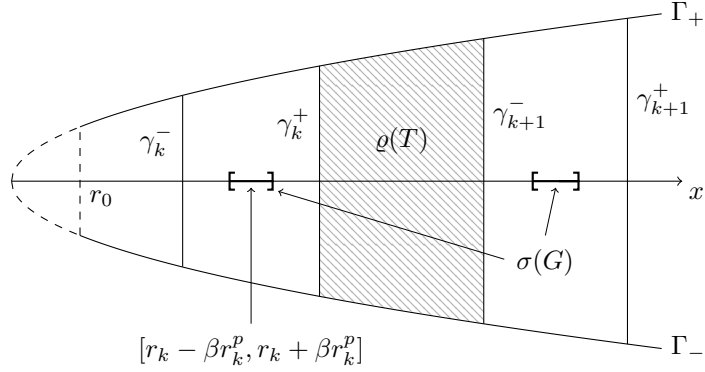
\begin{prop}\label{prop_specgap}
  Let\, $G$ be normal with compact resolvent,
  $\sigma(G)\cap\Omega(2\varphi)\subset\bbR_{\geq0}$ with 
  $0<\varphi\leq\pi/2$, $S$ $p$-subordinate
  to $G$ with bound $b$, $0\leq p<1$, and $T=G+S$. 

  Assume that there is
  a sequence $(r_k)_{k\geq 1}$ of nonnegative numbers tending monotonically
  to infinity and 
  constants $\beta\geq 0$, $\alpha>b$, $l>\beta+\alpha$ such that
  \begin{equation}\label{prop_specgap-specg}
    \sigma(G)\cap\bbR_{\geq0}\subset\bigcup_{k\geq1}[r_k-\beta r_k^p,
    r_k+\beta r_k^p]
  \end{equation}
  and
  \[r_k+lr_k^p \leq r_{k+1}-lr_{k+1}^p\]
  for almost all\, $k$. Then
  there are constants $C,r_0>0$, $k_0\geq 1$ such that the following holds:
  \begin{itemize}
  \item[(i)] The contours 
    $\Gamma_\pm$ from (\ref{lem_estim_contour1-contdef}) and
    \[\gamma_k^\pm=\bigl\{x+iy\,\big|\,x=r_k\pm lr_k^p,\,
      |y|\leq \alpha x^p\bigr\}\quad\text{with}\quad k\geq k_0\]
    as well as the regions enclosed by 
    $\gamma_k^+,\gamma_{k+1}^-,\Gamma_+,\Gamma_-$
    belong to $\varrho(T)$.
  \item[(ii)] $z\in\overline{\Omega(\varphi)}$ with $\Real z=r_k+lr_k^p$,
    $k\geq k_0$, implies $\|(T-z)^{-1}\|\leq C$.
  \item[(iii)] If\, $P_k$ and $Q_k$ are the Riesz projections of\,
    $T$ and $G$, respectively, associated with the region
    enclosed by $\gamma_k^-,\gamma_k^+,\Gamma_+,\Gamma_-$, then
    \[\sum_{k=k_0}^\infty|(P_ku|v)|\leq C\|u\|\|v\|\quad\text{for all}\quad
      u,v\in H\]
    and
    \[\dim\range(P_k)=\dim\range(Q_k)\quad\text{for}\quad k\geq k_0.\]
  \end{itemize}
\end{prop}
\textit{Proof.}
We set $s_k^\pm=r_k\pm lr_k^p$ so that
$r_k\leq s_k^+\leq s_{k+1}^-\leq r_{k+1}$.
Consider $s\in[s_k^+,s_{k+1}^-]$ with $k\geq k_0$. Then
\[s+\alpha s^p\leq s_{k+1}^-+\alpha r_{k+1}^p=r_{k+1}-(l-\alpha)r_{k+1}^p
  \leq r_{k+1}-\beta r_{k+1}^p.\]
Furthermore we have
\[s-\alpha s^p\geq s_k^+-\alpha (s_k^+)^p\]
for $k_0$ large enough, since the left-hand side is monotonically increasing
in $s$ for large $s$. In addition, the equivalent inequalities
\[
  s_k^+-\alpha (s_k^+)^p\geq r_k+\beta r_k^p \;\Leftrightarrow\;
  lr_k^p-\alpha(r_k+lr_k^p)^p\geq\beta r_k^p
\]
hold for $k_0$ sufficiently large since $1+lr_k^{p-1}\to 1$.
Using (\ref{prop_specgap-specg}), we have thus proved that,
for $k\geq k_0$,
\[s\in[s_k^+,s_{k+1}^-]\quad\Rightarrow\quad
  ]s-\alpha s^p,s+\alpha s^p[\,\subset\varrho(G).\]
With $r_0$ and $k_0$ appropriately chosen, 
Lemma~\ref{lem_psub_angle_resolvent} implies that the region
enclosed by $\gamma_k^+$, $\gamma_{k+1}^-$, $\Gamma^+$, and $\Gamma^-$
as well as the contours itself belong to $\varrho(T)$ for $k\geq k_0$.
Moreover, $\|(T-z)^{-1}\|$ and
$\|S(T-z)^{-1}\|$ are uniformly bounded for $z\in\overline{\Omega(\varphi)}$
with $\Real z=s_k^+$, $k\geq k_0$.
We also have $\dist(s_k^+,G)\geq\alpha(s_k^+)^p$ and
\[s_{k+1}^+-s_k^+=r_{k+1}-r_k+l(r_{k+1}^p-r_k^p)\geq 2lr_{k+1}^p.\]
The mean value theorem then yields
\[(s_{k+1}^+)^{1-p}-(s_k^+)^{1-p}\geq(1-p)(s_{k+1}^+)^{-p}(s_{k+1}^+-s_k^+)
  \geq \frac{2l(1-p)r_{k+1}^p}{\left(r_{k+1}+lr_{k+1}^p\right)^p},\]
i.e., $(s_{k+1}^+)^{1-p}-(s_k^+)^{1-p}\geq l(1-p)$ for $k\geq k_0$, $k_0$
sufficiently large.
We can thus apply Lemma~\ref{lem_rieszproj_estim} with $x_k=s_k^+$ to get 
the estimate for the sum over the Riesz projections.
The final claim is a consequence of Lemmas~\ref{lem_generalpert}
and \ref{lem_psub_angle_resolvent}.
\theoremend

\section{Existence of Riesz bases of invariant subspaces}
\label{sec_psubpert}

Let $G$ be an operator with compact resolvent.
Recall that we denote by $N_+(r_1,r_2,G)$ the sum of the multiplicities
of the eigenvalues of $G$ in the interval $]r_1,r_2[$, see
(\ref{eigmultsum_int}). Similarly, we write
\begin{equation}\label{eigmulsum}
  N(r,G)=\!\!\sum_{\lambda\in\sigma_p(G)\cap\overline{B_r(0)}}
  \!\!\dim \rsub(\lambda)
\end{equation}
for the sum of the multiplicities of all the eigenvalues $\lambda$ with
$|\lambda|\leq r$ and
\begin{equation}\label{eigmulsum_set}
  N(K,G)=\!\!\sum_{\lambda\in\sigma_p(G)\cap K}\!\!\dim \rsub(\lambda)
  \quad\text{for every set}\quad K\subset\bbC.
\end{equation}

Our first existence result for Riesz bases of invariant subspaces improves
a theorem due to 
Markus and Matsaev~\cite{markus-matsaev},
\cite[Theorem~6.12]{markus}; there, condition 
(\ref{theo_psubpert_invl2decomp-c}) was formulated with $\limsup$ 
instead of $\liminf$.
\begin{theo}
  \label{theo_psubpert_invl2decomp}
  Let\, $G$ be a normal operator with compact resolvent whose spectrum lies on
  a finite number of rays from the origin. Let\, $S$ be
  $p$-subordinate  to $G$ with $0\leq p<1$. If
  \begin{equation}\label{theo_psubpert_invl2decomp-c}
    \liminf_{r\to\infty}\frac{N(r,G)}{r^{1-p}}<\infty,
  \end{equation}
  then $T=G+S$ admits a Riesz basis of finite-dimensional
  $T$-invariant subspaces.
\end{theo}
\textit{Proof.}
Let $e^{i\theta_j}\bbR_{\geq 0}$ with $0\leq\theta_1<\ldots<\theta_n<2\pi$ be
the rays containing the eigenvalues of $G$ and let $S$ be $p$-subordinate
to $G$ with bound $b$.
From Theorem~\ref{theo_psubpert_specshape} we know that $T$ has a compact
resolvent and that almost all of its eigenvalues lie inside sectors of the 
form
\[\Omega_j=\bigl\{z\in\bbC\,\big|\,|\arg z-\theta_j|<\psi_j\bigr\}\quad
  \text{with}\quad 0<\psi_j\leq\frac{\pi}{4},\]
where the $\psi_j$ can be chosen such that these sectors are disjoint.
Lemma~\ref{lem_psub_angle_resolvent} shows that $\|(T-z)^{-1}\|$ is uniformly
bounded for $z\not\in\Omega_1\cup\ldots\cup\Omega_n$, $|z|\geq r_0$.
From the assumption (\ref{theo_psubpert_invl2decomp-c}) it can be
shown that
for each sector $\Omega_j$ there is a sequence $(r_{jk})_{k\geq 1}$ of 
positive numbers tending 
monotonically to infinity such that
\[\sup_kN_+(r_{jk}-2br_{jk}^p,r_{jk}+2br_{jk}^p,e^{-i\theta_j}G)<\infty,\]
see \cite[Lemma~6.11]{markus}.
By Proposition~\ref{prop_rieszproj_estim} we thus obtain a corresponding
sequence $(x_{jk})_{k\geq 1}$ such that $\|(T-z)^{-1}\|$ is 
uniformly bounded for
$z\in\Omega_j$, $\Real(e^{-i\theta_j}z)=x_{jk}$. 
Corollary~\ref{coroll_complete_rootspace} implies
that the system of root subspaces of $T$ is complete.

Furthermore, if $(P_{jk})_{k\geq 1}$ are the Riesz projections from
Proposition~\ref{prop_rieszproj_estim} corresponding to the eigenvalues
$\lambda\in\Omega_j$ of $T$ with $\Real(e^{-i\theta_j}\lambda)>x_{j1}$ and
$P_0$ is the Riesz projection for the (finitely many) remaining eigenvalues,
then
\[|(P_0u|v)|+\sum_{j=1}^n\sum_{k=1}^\infty|(P_{jk}u|v)|\leq C\|u\|\|v\|\]
with some constant $C\geq 0$.
Now Proposition~\ref{prop_proj_l2decomp} shows that the ranges of the 
projections $P_0$, $(P_{jk})_{j,k}$
form a Riesz basis.
\theoremendskip

Replacing condition \eqref{theo_psubpert_invl2decomp-c} by an assumption
on the localisation of the spectrum of $G$ on the rays,
we obtain our second perturbation theorem.
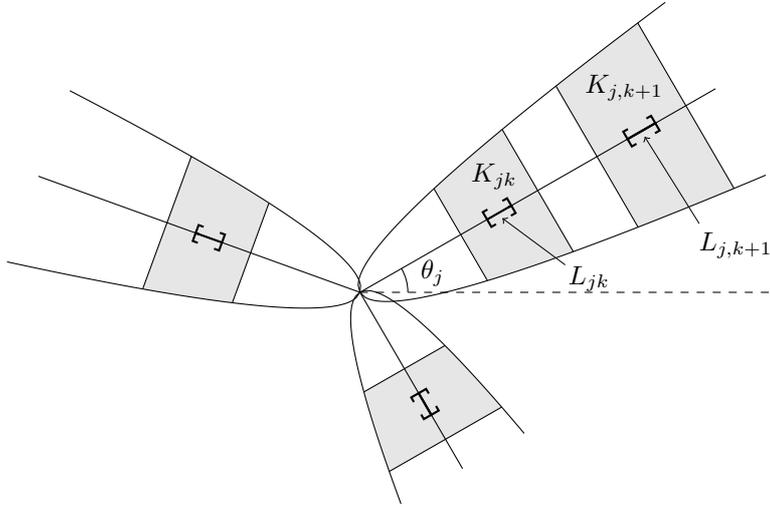
\begin{figure}
  \begin{center}
    \begin{tikzpicture}[scale=.9]
      \ifx\tikzdontfill\undefined
      \fill[rotate=-60,color=black!10]
        plot[id=theopsub_u11,samples=50,domain=1.3:2.5] function{.6*sqrt(x)}
	--(2.5,.95)--(2.5,-.95)
	--plot[id=theopsub_d11,samples=50,domain=1.3:2.5] function{-.6*sqrt(x)}
	--(1.3,-.68)--(1.3,.68)--cycle;
      \fill[rotate=160,color=black!10]
        plot[id=theopsub_u21,samples=50,domain=1.7:3] function{.6*sqrt(x)}
	--(3,1.04)--(3,-1.04)
	--plot[id=theopsub_d21,samples=50,domain=1.7:3] function{-.6*sqrt(x)}
	--(1.7,-.78)--(1.7,.78)--cycle;
      \fill[rotate=30,color=black!10]
        plot[id=theopsub_u31,samples=50,domain=1.7:3] function{.6*sqrt(x)}
	--(3,1.04)--(3,-1.04)
	--plot[id=theopsub_d31,samples=50,domain=1.7:3] function{-.6*sqrt(x)}
	--(1.7,-.78)--(1.7,.78)--cycle;
      \fill[rotate=30,color=black!10]
        plot[id=theopsub_u32,samples=50,domain=4:5.5] function{.6*sqrt(x)}
	--(5.5,1.41)--(5.5,-1.41)
	--plot[id=theopsub_d32,samples=50,domain=4:5.5] function{-.6*sqrt(x)}
	--(4,-1.2)--(4,1.2)--cycle;
      \fi

      \draw[dashed] (0,0)--(6,0);
      
      \draw (0,0)--(-60:3);
      \draw[rotate=-60] plot[id=theopsub_u1,samples=200,domain=0:3]
      function{.6*sqrt(x)};
      \draw[rotate=-60] plot[id=theopsub_d1,samples=200,domain=0:3] 
      function{-.6*sqrt(x)};
      
      \draw (0,0)--(160:5);
      \draw[rotate=160] plot[id=theopsub_u2,samples=200,domain=0:5]
      function{.6*sqrt(x)};
      \draw[rotate=160] plot[id=theopsub_d2,samples=200,domain=0:5] 
      function{-.6*sqrt(x)};
      
      \draw (0,0)--(30:6);
      \draw[rotate=30] plot[id=theopsub_u3,samples=200,domain=0:6]
      function{.6*sqrt(x)};
      \draw[rotate=30] plot[id=theopsub_d3,samples=200,domain=0:6] 
      function{-.6*sqrt(x)};
      
      \ifx\tikzdontfill\undefined
      \draw (.7,0) arc(0:30:.7);
      \draw (15:1.1) node{$\theta_j$};
      
      \draw[thick,[-] (-60:1.7)--(-60:2.1);
      \draw[thick,[-] (-60:2.1)--(-60:1.7);
      \draw[rotate=-60] (1.3,.68)--(1.3,-.68);
      \draw[rotate=-60] (2.5,.95)--(2.5,-.95);

      \draw[thick,[-] (160:2.13)--(160:2.57);
      \draw[thick,[-] (160:2.57)--(160:2.5);
      \draw[rotate=160] (1.7,.78)--(1.7,-.78);
      \draw[rotate=160] (3,1.04)--(3,-1.04);

      \draw[thick,[-] (30:2.13)--(30:2.57);
      \draw[thick,[-] (30:2.57)--(30:2.5);
      \draw[rotate=30] (1.7,.78)--(1.7,-.78);
      \draw[rotate=30] (3,1.04)--(3,-1.04);

      \draw[thick,[-] (30:4.5)--(30:5);
      \draw[thick,[-] (30:5)--(30:4.5);
      \draw[rotate=30] (4,1.2)--(4,-1.2);
      \draw[rotate=30] (5.5,1.41)--(5.5,-1.41);
      
      \draw (41:2.6) node{$K_{jk}$};
      \draw[rotate=30,very thin,->] (2.8,-1.15)--(2.35,-.1);
      \draw[rotate=30] (3,-1.5) node{$L_{jk}$};
      
      \draw (38:4.9) node{$K_{j,k+1}$};
      \draw[rotate=30,very thin,->] (4.8,-1.6)--(4.75,-.1);
      \fi
      \draw[rotate=30] (5.1,-2.15) node{$L_{j,k+1}$};

      
    \end{tikzpicture}
  \end{center}
  \caption{The situation of Theorem~\ref{theo_psubpert_rootl2decomp}}
  \label{figure_theopsub}
\end{figure}
\begin{theo}\label{theo_psubpert_rootl2decomp}
  Let\, $G$ be a normal operator with compact resolvent on a Hilbert space $H$
  and $S$ $p$-subordinate  to $G$ with bound $b$ and $0\leq p<1$.
  Suppose that the spectrum of\, $G$ lies on
  sequences of line segments on rays from the origin,
  \begin{equation}\label{theo_psubpert_rootl2decomp-sg}
    \sigma(G)\subset\bigcup_{j=1}^n\bigcup_{k\geq 1}L_{jk},\quad
    L_{jk}=\bigl\{e^{i\theta_j}x\,\big|\,x\geq 0,\,
    |x-r_{jk}|\leq\beta_j r_{jk}^p\bigr\},
  \end{equation}
  where $\beta_j\geq 0$, $0\leq\theta_1<\ldots<\theta_n<2\pi$, and
  $(r_{jk})_{k\geq 1}$ are monotonically increasing sequences of nonnegative
  numbers such that
  \begin{equation}\label{theo_psubpert_rootl2decomp-c}
    r_{jk}+l_jr_{jk}^p\leq r_{j,k+1}-l_jr_{j,k+1}^p
  \end{equation}
  for almost all\, $k$ with some constants $l_j>\beta_j+b$.

  Then $T=G+S$ has compact resolvent; 
  for $b<\alpha<\min\{l_1-\beta_1,\dots,l_n-\beta_n\}$
  almost all 
  eigenvalues of\, $T$ lie inside the regions 
  \begin{equation}\label{theo_psubpert_rootl2decomp-st}
    K_{jk}=\bigl\{e^{i\theta_j}(x+iy)\,\big|\,x\geq 0,\,|x-r_{jk}|\leq
    (\beta_j+\alpha)r_{jk}^p,\,|y|\leq \alpha x^p\bigr\},
  \end{equation}
  $j=1,\ldots,n$, $k\geq 1$ (cf.\ Figure~\ref{figure_theopsub});
  the spectral subspaces of\, $T$ corresponding to the regions $K_{jk}$
  together with the
  subspace corresponding to $\sigma(T)\setminus\bigcup_{j,k}K_{jk}$ form
  a Riesz basis of $H$; and we have
  \begin{equation}\label{theo_psubpert_rootl2decomp-em}
    N(L_{jk},G)=N(K_{jk},T)\quad\text{for almost all pairs }(j,k).
  \end{equation}

  Moreover, if there are constants $m,q>0$ such that for almost all
  pairs $(j,k)$ the assertions
  \begin{itemize}
    \item[(i)] $N(L_{jk},G)\leq m\quad$ and
    \item[(ii)] $\lambda_1,\lambda_2\in\sigma(T)\cap K_{jk},\,
      \lambda_1\neq\lambda_2
      \quad\Rightarrow\quad |\lambda_1-\lambda_2|>qr_{jk}^p$
  \end{itemize}
  hold, then the root subspaces of\, $T$ form a Riesz basis of\, $H$.
\end{theo}
\textit{Proof.}
We apply Theorem~\ref{theo_psubpert_specshape} and, for each ray,
Proposition~\ref{prop_specgap} with $\alpha$ and $l$ replaced by 
$\widetilde{\alpha}=(\alpha+b)/2$
and $\tilde{l}_j=\beta_j+\alpha$, respectively.
This shows that $T$ has a compact resolvent and that almost all of its
eigenvalues lie inside regions 
\[\bigl\{e^{i\theta_j}(x+iy)\,\big|\,x\geq 0,\,
  |x-r_{jk}|<\tilde{l}_jr_{jk}^p,\,
  |y|<\widetilde{\alpha} x^p\bigr\}\subset K_{jk}.\]
By Lemma~\ref{lem_psub_angle_resolvent}, $\|(T-z)^{-1}\|$ is
uniformly bounded outside certain disjoint sectors $\Omega_j$ around the rays
for $|z|$ large enough. 
For each ray, Proposition~\ref{prop_specgap} yields
a sequence $(x_{jk})_{k\in\bbN}$ tending monotonically to infinity such that
$\|(T-z)^{-1}\|$ is bounded for $z\in\Omega_j$, 
$\Real(e^{-i\theta_j}z)=x_{jk}$.
With Corollary~\ref{coroll_complete_rootspace} we conclude that
the system of root subspaces of $T$ is complete.
Moreover, we have
\[|(P_0u|v)|+\sum_{j=1}^n\sum_{k=1}^\infty|(P_{jk}u|v)|\leq C\|u\|\|v\|\]
for some $C\geq 0$ where $P_{jk}$ is the Riesz projection associated with
$K_{jk}$ and $P_0$ the one associated with
$\sigma(T)\setminus\bigcup_{jk}K_{jk}$;
Proposition~\ref{prop_proj_l2decomp} then yields the Riesz basis property.
Finally, if $Q_{jk}$ is the spectral projection of $G$ associated 
with $L_{jk}$, Proposition~\ref{prop_specgap} implies
$\dim\range(Q_{jk})=\dim\range(P_{jk})$ for almost all $(j,k)$ and thus
(\ref{theo_psubpert_rootl2decomp-em}).

Now suppose that with $m,q>0$ the additional assumptions (i) and (ii) hold
for almost all pairs $(j,k)$.
We aim to show that the root subspaces corresponding to the eigenvalues
of $T$ in each $K_{jk}$ form a Riesz basis of $\range(P_{jk})$
with constant $c$ independent of $(j,k)$.
Without loss of generality we may assume
\[\theta_j=0 \quad\text{and}\quad
  \beta_j+\alpha+q\leq l_j.\]
We want to apply Corollary~\ref{coroll_resolvent_discs2}
with $l_0=\beta_j+\alpha$ and set $l$ accordingly. 
From the relation $r_{j,k+1}-r_{jk}\geq 2l_jr_{jk}^p$ it is easy to verify
that the number of elements
$r_{jk}$ in the interval $[r-lr^p,r+lr^p]$ is
at most $2l/l_j$ for $r$ sufficiently large. Hence there is a
constant $m_0$ such that 
\[N_+(r-lr^p,r+lr^p,G)\leq m_0\quad\text{for $r$ sufficiently large.}\]
Let $\lambda$ be an eigenvalue of $T$ in $K_{jk}$. By 
Corollary~\ref{coroll_resolvent_discs2} there exists $q_1\in\,]0,q[$
such that the points $w$ on the circle around $\lambda$ with radius
$q_1r_{jk}^p$ satisfy $\|(T-w)^{-1}\|\leq C_0C_1^{m_0}r_{jk}^{-p}$.
In addition, this circle lies inside the strip
$|\Real z-r_{jk}|\leq l_jr_{jk}^p$ and assumption (ii) thus 
implies that $\lambda$ is the only possible eigenvalue of $T$ inside
that circle. 
Therefore, the Riesz projection $P_\lambda$ for $\lambda$
satisfies
\[\|P_\lambda\|\leq2\pi q_1r_{jk}^p\frac{C_0C_1^{m_0}}{r_{jk}^p}
  \leq2\pi qC_0C_1^{m_0}.\]
If $\lambda_1,\ldots,\lambda_{m_1}$ are the eigenvalues of $T$
in $K_{jk}$, we have $m_1\leq N(K_{jk},T)\leq m$ and conclude
\[\sum_{s=1}^{m_1}|(P_{\lambda_s}u|v)|\leq2\pi mqC_0C_1^{m_0}\|u\|\|v\|.\]
According to Proposition~\ref{prop_proj_l2decomp}, the subspaces
$\range(P_{\lambda_s})$, $s=1,\ldots,m_1$, form a Riesz basis of 
$\range(P_{jk})$ with constant $c$ independent of $k$. This is true
for almost all pairs $(j,k)$, and hence an application of 
Lemma~\ref{lem_join_l2decomp}
shows that the root subspaces of $T$ form a Riesz basis of $H$.
\theoremend

\begin{remark}\label{rem_psubpert_basiseigvec}
  If almost all eigenvalues of $G$ are simple and almost all
  line segments $L_{jk}$ contain one eigenvalue only, then 
  Theorem~\ref{theo_psubpert_rootl2decomp} yields a Riesz basis of
  eigenvectors and finitely many Jordan chains for $T$.
  Indeed almost all spectral subspaces corresponding to the $K_{jk}$
  are one-dimensional in this case, and the
  Riesz basis of subspaces is thus equivalent to the existence of a
  Riesz basis of eigenvectors and finitely many Jordan chains.
  \definend
\end{remark}

\begin{remark}\label{rem_spec_pasymp2density}
  It can be shown \cite[Lemma~3.4.9]{wyss} 
  that, if $G$ satisfies the spectral condition
  (\ref{theo_psubpert_rootl2decomp-sg}) with some $\beta_j>0$ such that
  $r_{jk}+\beta_jr_{jk}^p\leq r_{j,k+1}-\beta_jr_{j,k+1}^p$
  and $N(L_{jk},G)$ is bounded in $(j,k)$, then 
  $\sup_{r\geq 1}N(r,G)r^{p-1}<\infty$; in particular
  the spectral condition (\ref{theo_psubpert_invl2decomp-c}) 
  of Theorem~\ref{theo_psubpert_invl2decomp}  holds.
  However, the first part of Theorem~\ref{theo_psubpert_rootl2decomp} 
  is applicable even if $N(L_{jk},G)$ is unbounded and 
  (\ref{theo_psubpert_invl2decomp-c}) does not hold.
  \definend
\end{remark}

The condition (\ref{theo_psubpert_rootl2decomp-c}) can be reformulated
for sequences with a certain asymptotic behaviour:
\begin{lemma}\label{lem_asymp2gap}
  Consider the sequence of nonnegative numbers given by 
  \[r_k=ck^q+d_kk^{q-1}\]
  with $c>0$, $q\geq 1$ and a converging sequence $(d_k)_{k\in\bbN}$. 
  Then for $l,p\geq 0$ the relation
  \[r_k+lr_k^p\leq r_{k+1}-lr_{k+1}^p\]
  holds for almost all\, $k\in\bbN$ if
  \begin{itemize}
  \item[(i)] $p<1-1/q$, or
  \item[(ii)] $p=1-1/q$ and $l<qc^{1/q}/2$.
  \end{itemize}
\end{lemma}
\textit{Proof.}
This can be shown in a straightforward way by a Taylor series expansion
of $r_{k+1}$ in $k$.
\theoremendskip

The next proposition reverses the assertions of the
Theorems~\ref{theo_psubpert_specshape},
\ref{theo_psubpert_invl2decomp} and
\ref{theo_psubpert_rootl2decomp} to some extend.
As a consequence, the assumptions in these theorems can be relaxed.
\begin{prop}\label{prop_psubpert_revert}
  Let\, $G$ be an operator on a Hilbert space $H$
  with compact resolvent and a Riesz basis
  of Jordan chains. Suppose that\, $0\leq p<1$, $\alpha\geq 0$,
  $0\leq\theta_j<2\pi$, $j=1,\ldots,n$, such that either
  \begin{itemize}
  \item[(i)] there exists $r_0>0$ with
    \[\sigma(G)\subset B_{r_0}(0)\cup\bigcup_{j=1}^n\bigl\{
      e^{i\theta_j}(x+iy)\,\big|\,x>0,\,|y|\leq\alpha x^p\bigr\},
      \quad\text{or}\]
  \item[(ii)] almost all eigenvalues of\, $G$ lie inside regions
    \[K_{jk}=\bigl\{e^{i\theta_j}(x+iy)\,\big|\,r_{jk}^-\leq x\leq r_{jk}^+,\,
      |y|\leq \alpha x^p\bigr\},\quad j=1,\ldots,n,\, k\geq 1,\]
    where $(r_{jk}^\pm)_{k\geq 1}$ are
    sequences of positive numbers satisfying
    \(r_{jk}^-\leq r_{jk}^+<r_{j,k+1}^-.\)
  \end{itemize}
  Then there is an isomorphism $J:H\to H$, a normal operator $G_0$ on $H$
  with compact resolvent, and an operator $S_0$ 
  $p$-subordinate  to $G_0$ such that
  \[J\mdef(G)=\mdef(G_0),\quad JGJ^{-1}=G_0+S_0.\]
  In case (i), all eigenvalues of\, $G_0$ lie on the rays 
  $e^{i\theta_j}\bbR_{\geq 0}$ and we have
  \[N(r,G_0)=N(r,G)\quad\text{for}\quad r\geq 1.\]
  In case (ii), all eigenvalues of\, $G_0$ lie on the line segments
  \[L_{jk}=\bigl\{e^{i\theta_j}x\,\big|\,r_{jk}^-\leq x\leq r_{jk}^+\bigr\},\]
  and $N(L_{jk},G_0)=N(K_{jk},G)$ holds for almost all pairs $(j,k)$.
  
  Moreover, if\, $S$ is $p$-subordinate  to $G$, then $JSJ^{-1}$
  is $p$-sub\-or\-di\-nate  to $G_0$.
\end{prop}
\textit{Proof.}
  The idea is to use the isomorphism $J$ to transform
  the Riesz basis of Jordan chains of $G$ to an orthonormal basis
  of eigenvectors of $G_0$.
  Then one associates with each eigenvalue $\lambda=e^{i\theta_j}(x+iy)$ of $G$
  an eigenvalue $\mu=e^{i\theta_j}w$ of $G_0$ with $w>0$ such that the 
  assertions on the spectrum hold.
A complete proof can be found in \cite[\S3.4]{wyss}.
\theoremend

\begin{remark}\label{rem_psubpert_revert}
  Theorems~\ref{theo_psubpert_specshape} and
  \ref{theo_psubpert_invl2decomp} also hold if $G$ is as in the
  previous proposition and satisfies condition 
  \ref{prop_psubpert_revert}(i).
  Indeed we have
  \[J(G+S)J^{-1}=G_0+S_0+JSJ^{-1}\]
  in this case, $S_0+JSJ^{-1}$ is $p$-subordinate to $G_0$,
  and the theorems can be applied to the right-hand side.
  Analogously, Theorem~\ref{theo_psubpert_rootl2decomp}  also holds
  if $G$ satisfies \ref{prop_psubpert_revert}(ii).
  In both cases, $b$ is now the $p$-subordination bound of $S_0+JSJ^{-1}$
  to $G_0$.
\end{remark}

\section{Application to block operator matrices}
\label{sec_psubpert_appli}

We apply Theorems~\ref{theo_psubpert_invl2decomp} and 
\ref{theo_psubpert_rootl2decomp} to two classes of
diagonally dominant block operator matrices.
For many results about the spectral theory of block operator matrices
see \cite{tretter2,tretter-book}.
\begin{theo}\label{theo_psubblkop}
  Let\, $A(H_1\to H_1)$ and $D(H_2\to H_2)$ be normal operators with
  compact resolvents on Hilbert spaces such that 
  the spectra of\, $A$ and $D$ lie on finitely many rays from the origin
  and
  \[\liminf_{r\to\infty}\frac{N(r,A)}{r^{1-p}}<\infty,\quad
    \liminf_{r\to\infty}\frac{N(r,D)}{r^{1-p}}<\infty\]
  with $0\leq p<1$.
  Suppose that the operators $C(H_1\to H_2)$ and $B(H_2\to H_1)$ 
  are $p$-subordinate\footnote{
    This notion of $p$-subordination is more general than the one from
    Section~\ref{sec_psubspec}, since the operators $B$ and $C$ map
    from one Hilbert space into a (possibly) different one.}
  to $A$ and $D$, respectively,
  \begin{align*}
    &\|Cu\|\leq b\|u\|^{1-p}\|Au\|^p \quad\text{for}\quad
      u\in\mdef(A)\subset\mdef(C),\\
    &\|Bv\|\leq b\|v\|^{1-p}\|Dv\|^p \quad\text{for}\quad
      v\in\mdef(D)\subset\mdef(B).
  \end{align*}
  Then the block operator matrix
  \[T=\pmat{A&B\\C&D}\]
  acting on $H_1\times H_2$
  has a compact resolvent, admits a Riesz basis of 
  finite-dimen\-sional $T$-invariant subspaces,
  and for every $\alpha>b$ there is a constant\, $r_0\geq 0$ such that
  \[\sigma(T)\subset B_{r_0}(0)\cup\bigcup_{j=1}^{n}\bigl\{e^{i\theta_j}(x+iy)
    \,\big|\,x\geq 0,|y|\leq\alpha x^p\bigr\}.\]
  Here $\theta_1,\ldots,\theta_n$ with $0\leq\theta_j<2\pi$ are the angles of 
  the rays on which the spectra of\, $A$ and $D$ lie.
\end{theo}
\textit{Proof.}
We apply Theorems~\ref{theo_psubpert_specshape} and
\ref{theo_psubpert_invl2decomp} to the decomposition
\[T=G+S\quad\text{with}\quad G=\pmat{A&0\\0&D},\quad S=\pmat{0&B\\C&0}.\]
Indeed $G$ is normal with compact resolvent and
\[\sigma(G)=\sigma(A)\cup\sigma(D), \quad
  N(r,G)=N(r,A)+N(r,D).\]
Moreover, using H\"older's inequality, we find
\begin{align*}
  \Big\|S\pmat{u\\v}\Big\|^2&=\|Bv\|^2+\|Cu\|^2\leq b^2\|v\|^{2(1-p)}
  \|Dv\|^{2p}+b^2\|u\|^{2(1-p)}\|Au\|^{2p}\\
  &\leq b^2\left(\|u\|^2+\|v\|^2\right)^{1-p}\left(\|Au\|^2
  +\|Dv\|^2\right)^p  
\end{align*}
for $u\in\mdef(A)$, $v\in\mdef(D)$, i.e.
\[\|Sw\|\leq b\|w\|^{1-p}\|Gw\|^p \quad\text{for}\quad
  w\in\mdef(G)=\mdef(A)\times\mdef(D);\]
$S$ is $p$-subordinate to $G$.
\theoremendskip

In the next theorem, a symmetry of the operator matrix with respect
to an indefinite inner product yields a gap in the spectrum
around the imaginary axis.
This makes it possible to
apply the second part of Theorem~\ref{theo_psubpert_rootl2decomp}.
\begin{theo}\label{theo_uposblkop}
  Let $A$ be a skew-adjoint operator with compact resolvent on a 
  Hilbert space $H$.
  Let $B,C:H\to H$ be bounded, selfadjoint and uniformly positive,
  $B,C\geq\gamma>0$. 
  Write $(ir_k)_{k\in\Lambda}$ for the sequence of eigenvalues of $A$
  where $\Lambda\in\{\bbZ_+,\bbZ_-,\bbZ\}$ 
  and $(r_k)_k$ is increasing.
  Suppose that almost all eigenvalues $ir_k$ are simple and that
  for some $l>b=\max\{\|B\|,\|C\|\}$ we have
  \[r_{k+1}-r_k\geq 2l \quad\text{for almost all\, }k\in\Lambda.\]
  Then the block operator matrix
  \[T=\pmat{A&B\\C&A}\]
  has a compact resolvent, its spectrum is symmetric with respect
  to the imaginary axis and satisfies
  \[\sigma(T)\subset\bigl\{z\in\bbC\,\big|\,|z-ir_k|\leq b
    \text{ for some }k,\,|\Real z|\geq\gamma\bigr\}.\]
  Moreover almost all eigenvalues are simple and
  $T$ admits a Riesz basis of eigenvectors and finitely many
  Jordan chains.
\end{theo}
\textit{Proof.}
We consider the decomposition
\[T=G+S,\quad G=\pmat{A&0\\0&A},\quad S=\pmat{0&B\\C&0}.\]
$G$ is skew-adjoint with compact resolvent, 
$\sigma(G)=\{ir_k\,|\,k\in\Lambda\}$, and almost all of its 
eigenvalues have multiplicity $2$. $S$ is bounded with 
$\|S\|=b$. 
By Theorem~\ref{theo_psubpert_specshape}
$T$ has a compact resolvent. If $z$ is a point outside the discs $D_k$
with radius
$b$ around the $ir_k$, then $\dist(z,\sigma(G))>b$ and $\|(G-z)^{-1}\|<b^{-1}$;
thus $z\in\varrho(T)$ by Lemma~\ref{lem_generalpert}.

Now we use the indefinite inner products $(J_j\cdot|\cdot)$ on $H\times H$
given by the fundamental symmetries
\[J_1=\pmat{0&-iI\\iI&0},\qquad J_2=\pmat{0&I\\I&0},\]
where $(\cdot|\cdot)$ is the standard scalar product on $H\times H$.
We refer to \cite{azizov-iokhvidov} for a treatment of
indefinite inner product spaces and operators therein.
It is easy to verify that $T$ is $J_1$-skew-adjoint
(i.e., $J_1T$ is skew-adjoint), which implies that
$\sigma(T)$ is symmetric with respect to $i\bbR$.
On the other hand, for an eigenvalue $\lambda$ of $T$ with
eigenvector $w$ an easy calculation yields
\[\gamma\|w\|^2\leq\Real(J_2Tw|w)\leq|\Real\lambda|
  |(J_2w|w)|\leq|\Real\lambda|\|w\|^2;\]
hence $|\Real\lambda|\geq\gamma$, which shows the asserted shape 
of the spectrum.

Finally we apply Theorem~\ref{theo_psubpert_rootl2decomp} with 
$p=\beta_1=\beta_2=0$, $\theta_1=\pi/2$, $\theta_2=3\pi/2$.
It shows that $N(D_k,T)=2$ for almost all discs $D_k$.
Consequently, almost all $D_k$ contain only one skew-conjugate pair 
of simple eigenvalues
$\lambda,\,-\overline{\lambda}$ with $|\Real\lambda|\geq\gamma$.
The second part of the theorem thus
implies that the root subspaces of $T$ form a
Riesz basis. Since almost all root subspaces have dimension
one, this is equivalent to the existence of a
Riesz basis of eigenvectors and finitely many Jordan chains.
\theoremend

\section*{Acknowledgement}
The author wishes to thank Alexander Markus for some comments on
the conditions in Theorem~\ref{theo_psubpert_invl2decomp}
and Christiane Tretter for many valuable suggestions.
{The author is also grateful for the support of
    Deutsche Forschungsgemeinschaft DFG, grant no.\ TR368/6-1, and
    of Schweizerischer Nationalfonds SNF, grant no.~15-486.}

\bibliographystyle{cwyss}
\bibliography{biblist}

\end{document}